\documentclass[11pt]{article}

\usepackage{amssymb, amsfonts, amsmath, amsthm, graphics} 

\newtheorem{theorem}{Theorem}[section]
\newtheorem{lemma}[theorem]{Lemma}
\newtheorem{proposition}[theorem]{Proposition}
\newtheorem{corollary}[theorem]{Corollary}

\theoremstyle{definition}
\newtheorem{definition}[theorem]{Definition}
\newtheorem{notation}[theorem]{Notation}
\newtheorem{remark}[theorem]{Remark}
\newtheorem{example}[theorem]{Example}

\newcommand{\cA}{ {\mathcal A} }
\newcommand{\cAo}{ {\mathcal A}^o }
\newcommand{\cB}{ {\mathcal B} }
\newcommand{\cF}{ {\mathcal F} }
\newcommand{\cM}{ {\mathcal M} }

\newcommand{\cS}{ {\mathcal S} }
\newcommand{\cT}{ {\mathcal T} }
\newcommand{\cV}{ {\mathcal V} }
\newcommand{\cW}{ {\mathcal W} }
\newcommand{\cX}{ {\mathcal X} }
\newcommand{\cY}{ {\mathcal Y} }

\newcommand{\bC}{ {\mathbb C} }

\newcommand{\bG}{ {\mathbb G} }

\newcommand{\bR}{ {\mathbb R} }
\newcommand{\bZ}{ {\mathbb Z} }

\newcommand{\fD}{ {\mathfrak D} }
\newcommand{\fM}{ {\mathfrak M} }

\newcommand{\abs}{ \mbox{Abs} }
\newcommand{\mult}{ \mbox{Mult} }

\newcommand{\Kr}{ \mbox{Kr} }
\newcommand{\bo}{ \mbox{Bo} }
\newcommand{\so}{ \mbox{So} }
\newcommand{\moeb}{ \mbox{M\"ob} }
\newcommand{\moebA}{ {\mbox{M\"ob}}^{(A)} }
\newcommand{\moebB}{ {\mbox{M\"ob}}^{(B)} }
\newcommand{\ncb}{ NC^{(B)} }
\newcommand{\nczb}{ NCZ^{(B)} }
\newcommand{\Dno}{ D_{\#} }
\newcommand{\muno}{ \mu_{\#} }

\newcommand{\ee}{\varepsilon}
\newcommand{\kk}{ \kappa }
\newcommand{\ckn}{ \underline{\kappa}_{n} }
\newcommand{\kappaA}{ \kappa^{(A)} }
\newcommand{\kappaB}{ \kappa^{(B)} }
\newcommand{\ktild}{ \widetilde{\kappa} }
\newcommand{\cktild}{ \underline{\widetilde{\kappa}}_{n} }
\newcommand{\phiA}{ \varphi^{(A)} }
\newcommand{\phiB}{ \varphi^{(B)} }
\newcommand{\phiprim}{\varphi '}
\newcommand{\psiprim}{\psi '}
\newcommand{\alphatild}{ \widetilde{\alpha} }
\newcommand{\phitild}{ \widetilde{\varphi} }
\newcommand{\psitild}{ \widetilde{\psi} }

\newcommand{\piciup}{ \stackrel{\vee}{\pi} }
\newcommand{\freestar}{ \framebox[7pt]{$\star$} }

\begin{document}

$\ $

\begin{center}
{\bf\Large Infinitesimal non-crossing cumulants}

\vspace{6pt}

{\bf\Large and free probability of type B}

\vspace{20pt}

{\large Maxime F\'evrier \hspace{2cm}
Alexandru Nica \footnote{Research supported by a Discovery Grant
from NSERC, Canada.} }

\vspace{10pt}

\end{center}

\begin{abstract}
Free probabilistic considerations of type B first appeared 
in the paper of Biane, Goodman and Nica [Trans. AMS 355 (2003), 
2263-2303]. Recently, connections between type B and 
infinitesimal free probability were put into evidence by Belinschi 
and Shlyakhtenko [arXiv:0903.2721]. The interplay between 
``type B'' and ``infinitesimal'' is also the object of the 
present paper. We study infinitesimal freeness for a family
of unital subalgebras $\cA_1, \ldots , \cA_k$ in an 
infinitesimal noncommutative probability space
$( \cA , \varphi , \phiprim )$ and we introduce a concept of 
{\em infinitesimal non-crossing cumulant functionals} for 
$( \cA , \varphi , \phiprim )$, obtained by taking a formal 
derivative in the formula for usual non-crossing cumulants.
We prove that the infinitesimal freeness of $\cA_1, \ldots , \cA_k$
is equivalent to a vanishing condition for mixed cumulants; this 
gives the infinitesimal counterpart for a theorem 
of Speicher from ``usual'' free probability. 
We show that the lattices $\ncb (n)$ of non-crossing partitions
of type B appear in the combinatorial study of  
$( \cA , \varphi , \phiprim )$, in the formulas for infinitesimal 
cumulants and when describing alternating products of 
infinitesimally free random variables. As an application of 
alternating free products, we observe the infinitesimal analogue
for the well-known fact that freeness is preserved under 
compression with a free projection. As another application, we 
observe the infinitesimal analogue for a well-known procedure 
used to construct free families of free Poisson elements.
Finally, we discuss situations when the freeness of 
$\cA_1, \ldots , \cA_k$ in $( \cA , \varphi )$ can be naturally 
upgraded to infinitesimal freeness in $( \cA , \varphi , \phiprim )$, 
for a suitable choice of a ``companion functional''
$\phiprim : \cA \to \bC$.
\end{abstract}

\vspace{6pt}

\begin{center}
{\bf\large 1. Introduction}
\setcounter{section}{1}
\end{center}

\begin{center}
{\bf 1.1 The framework of the paper}
\end{center}

This paper is concerned with a form of free independence for 
noncommutative random variables, which can be called 
``freeness of type B'' or ``infinitesimal freeness'', and 
occurs in relation to objects of the form
\begin{equation}   \label{eqn:1.11}
\left\{  \begin{array}{l}
( \cA , \varphi , \phiprim ), \ \mbox{ where $\cA$ is a 
                unital algebra over $\bC$}                   \\
\mbox{$\ \ \ $ and $\varphi , \phiprim : \cA \to \bC$ are 
             linear with $\varphi (1_{\cA}) = 1$,
             $\phiprim (1_{\cA} ) = 0$.                  }
\end{array}  \right.
\end{equation}
The motivation for considering objects as in (\ref{eqn:1.11}) 
is three-fold.

(a) This framework generalizes the link-algebra associated to 
a noncommutative probability space of type B, in the sense 
introduced by Biane, Goodman and Nica \cite{BGN2003}. One can
thus take the point of view that (\ref{eqn:1.11}) provides us
with an enlarged framework for doing ``free probability of 
type B''. This point of view is justified by the fact that 
lattices of non-crossing partitions of type B do indeed appear 
in the underlying combinatorics -- see e.g. Theorem \ref{thm:6.4}
below, concerning alternating products of infinitesimally free 
random variables.

(b) It turns out to be beneficial to consolidate the 
functionals $\varphi , \phiprim$ from (\ref{eqn:1.11}) into only 
one functional 
\begin{equation}   \label{eqn:1.12}
\phitild : \cA \to \bG, \ \ 
\phitild := \varphi + \ee \phiprim ,
\end{equation}
where $\bG$ denotes the two-dimensional Grassman algebra generated 
by an element $\ee$ which satisfies $\ee^2 =0$. Thus $\bG$ is the 
extension of $\bC$ defined as
\begin{equation}   \label{eqn:1.13}
\bG = \{ \alpha + \ee \beta \mid \alpha , \beta \in \bC \},
\end{equation}
with multiplication given by 
$( \alpha_1 + \ee \beta_1) \cdot ( \alpha_2 + \ee \beta_2 )$ = 
$\alpha_1 \alpha_2 + \ee ( \alpha_1 \beta_2 + \beta_1 \alpha_2)$, 
and the structure from (\ref{eqn:1.11}) could equivalently be 
treated as 
\begin{equation}   \label{eqn:1.14}
\left\{  \begin{array}{l}
( \cA , \phitild ), \ \mbox{ where $\cA$ is a 
                unital algebra over $\bC$}                   \\
\mbox{$\ \ \ $ and $\phitild : \cA \to \bG$ is $\bC$-linear
      with $\phitild (1_{\cA}) = 1$.}
\end{array}  \right.
\end{equation}
The framework (\ref{eqn:1.14}) was discussed in the PhD Thesis 
of Oancea \cite{O2007}, under the name of ``scarce 
\footnote{The adjective ``scarce'' is used in order to distinguish
from the concept of ``$\bG$-probability space'' from operator-valued 
free probability, where one would require the functional 
$\phitild$ to be $\bG$-linear.}
$\bG$-probability space''. Specifically, Chapter 7 of \cite{O2007}
considers a concept of $\bG$-freeness for a family of unital 
subalgebras in a $\bG$-probability space, which is defined 
via a vanishing condition for mixed $\bG$-valued cumulants, and 
generalizes the concept of freeness of type B from \cite{BGN2003}.

(c) The recent paper \cite{BS2009} by Belinschi and Shlyakhtenko 
discusses a concept of ``infinitesimal distribution''
$( \bC \langle X_1, \ldots , X_k \rangle , \mu , \mu ' )$ 
which is exactly as in (\ref{eqn:1.11}), with 
$\bC \langle X_1, \ldots , X_k \rangle$ denoting the algebra
of polynomials in noncommuting indeterminates 
$X_1, \ldots , X_k$. This remarkable paper brings forth the idea 
that interesting infinitesimal distributions arise when $\mu$ is 
the limit at $0$ and $\mu '$ is the derivative at $0$ for a 
family of $k$-variable distributions 
$( \mu_t : \bC \langle X_1, \ldots ,X_k \rangle \to \bC )_{t \in T}$,
where $T$ is a set of real numbers having 0 as accumulation point.
As we will show below, this ties in really nicely with the 
$\bG$-valued cumulant considerations mentioned in (b); indeed, one 
could say that \cite{BS2009} puts the $\ee$ from (\ref{eqn:1.13}) 
in its right place -- it is a sibling of the $\ee$'s from calculus, 
only that instead of just having ``$\ee^2$ much smaller than $\ee$'' 
one goes for the radical requirement that $\ee^2 = 0$.

Upon consideration, it seems that what goes best with the 
framework from (\ref{eqn:1.11}) is the ``infinitesimal'' terminology
from (c), which is in particular adopted in the next definition.
Throughout the paper some terminology inspired from (a) and (b) will 
also be used, in the places where it is suggestive to do so (e.g. 
when talking about ``soul companions for $\varphi$'' in subsection 
1.3 below).

\begin{definition}      \label{def:1.1}

$1^o$ A structure $( \cA , \varphi , \phiprim )$ as in 
(\ref{eqn:1.11}) will be called an {\em infinitesimal noncommutative 
probability space} (abbreviated as {\em incps}). 

$2^o$ Let $( \cA , \varphi , \phiprim )$ be an incps and let 
$\cA_1, \ldots , \cA_k$ be unital subalgebras of $\cA$. We will 
say that $\cA_1, \ldots , \cA_k$ are {\em infinitesimally free} 
with respect to $( \varphi , \phiprim )$ when they satisfy the 
following condition:

\begin{center}
If $i_1, \ldots , i_n \in \{ 1, \ldots , k \}$ are such that 
$i_1 \neq i_2, i_2 \neq i_3 , \ldots , i_{n-1} \neq i_n$, 

and if $a_1 \in \cA_{i_1}, \ldots , a_n \in \cA_{i_n}$ 
are such that $\varphi (a_1) = \cdots = \varphi (a_n) = 0$, 

then $\varphi ( a_1 \cdots a_n ) = 0$ and
\end{center}
\begin{equation}   \label{eqn:1.15}
\phiprim ( a_1 \cdots a_n ) = 
\left\{  \begin{array}{l}
\varphi (a_1 \, a_n) \varphi (a_2 \, a_{n-1}) \cdots 
       \varphi ( a_{(n-1)/2} \, a_{(n+3)/2} ) \cdot
       \phiprim ( a_{(n+1)/2} ),                                \\
\mbox{$\ \ $} \ \ \mbox{ if $n$ is odd and $i_1 = i_n , 
      i_2 = i_{n-1}, \ldots , i_{(n-1)/2} = i_{(n+3)/2}$, }     \\
0, \ \mbox{ otherwise.}
\end{array}  \right.
\end{equation}
\end{definition}

Recall that in the free probability literature it is customary to
use the name {\em noncommutative probability space} for a pair 
$( \cA , \varphi )$ where $\cA$ is a unital algebra over $\bC$ 
and $\varphi : \cA \to \bC$ is linear with $\varphi (1_{\cA}) =1$.
Thus the concept of {\em infinitesimal} noncommutative probability 
space is obtained by adding to $( \cA , \varphi )$ another 
functional $\phiprim$ as in (\ref{eqn:1.11}). It is also 
immediate that Definition \ref{def:1.1}.$2^o$ of infinitesimal 
freeness is obtained by adding the condition (\ref{eqn:1.15}) to
the ``usual'' definition for the freeness of $\cA_1, \ldots , \cA_k$
in $( \cA , \varphi )$ (as appearing e.g. in \cite{VDN1992}, 
Definition 2.5.1).

Definition \ref{def:1.1}.$2^o$ is a reformulation of the concept
with the same name from Definition 13 of \cite{BS2009}. The relations
with \cite{BS2009}, \cite{BGN2003} are discussed more precisely 
in Section 2 (cf. Remarks \ref{rem:2.8}, \ref{rem:2.9}). 
Section 2 also collects a few miscellaneous properties of 
infinitesimal freeness that follow directly from the definition. 
Most notable among them is that one can easily extend to 
infinitesimal framework the well-known free 
product construction of noncommutative probability spaces 
$( \cA_1, \varphi_1 ) * \cdots * ( \cA_k , \varphi_k )$, as 
presented e.g. in Lecture 6 of \cite{NS2006}. More precisely: if
$( \cA_1, \varphi_1 ) * \cdots * ( \cA_k, \varphi_k )$
$=: ( \cA , \varphi )$ and if we are given linear functionals 
$\phiprim_i : \cA_i \to \bC$ such that $\phiprim_i ( 1_{\cA} ) = 0$, 
$1 \leq i \leq k$, then there exists a unique linear functional 
$\phiprim : \cA \to \bC$ such that $\phiprim \mid \cA_i = \phiprim_i$,
$1 \leq i \leq k$, and such that $\cA_1, \ldots , \cA_k$ are 
infinitesimally free in $( \cA , \varphi , \phiprim )$. (See 
Proposition \ref{prop:2.4} below.) The resulting incps 
$( \cA , \varphi , \phiprim )$ can thus be taken, by definition,
as the {\em free product} of $( \cA , \varphi_i , \phiprim_i )$ 
for $1 \leq i \leq k$.

$\ $

\begin{center}
{\bf 1.2 Non-crossing cumulants for 
\boldmath{$( \cA , \varphi , \phiprim )$} }
\end{center}

An important tool in the combinatorics of free probability is 
the family of non-crossing cumulant functionals 
$( \kappa_n : \cA^n \to \bC )_{n \geq 1}$ associated to a 
noncommutative probability space $( \cA , \varphi )$. These 
functionals were introduced in \cite{S1994}; for a detailed 
presentation of their basic properties, see Lecture 11 of 
\cite{NS2006}. For every $n \geq 1$, the multilinear functional 
$\kappa_n : \cA^n \to \bC$ is defined via a certain summation 
formula over the lattice $NC(n)$ of non-crossing partitions of 
$\{ 1, \ldots , n \}$. We will review the formula for a general
$\kappa_n$ in subsection 3.2, here we only pick a special value of 
$n$ that we use for illustration, e.g. $n=3$. In this special 
case one has
\begin{equation}    \label{eqn:1.21}
\kappa_3 (a_1, a_2, a_3) = 
\begin{array}[t]{lr}
\varphi (a_1 a_2 a_3) - \varphi (a_1) \varphi (a_2 a_3)  
               - \varphi (a_2) \varphi (a_1 a_3)   &        \\
- \varphi (a_3) \varphi (a_1 a_2)
+ 2 \varphi (a_1) \varphi (a_2) \varphi (a_3),     & 
\ \ \forall \, a_1, a_2, a_3 \in \cA .
\end{array}
\end{equation}  
The expression on the right-hand side of (\ref{eqn:1.21}) has 
5 terms (premultiplied by integer coefficients 
\footnote{The meaning of these coefficients is that they 
are special values of the M\"obius function of $NC(3)$, as reviewed 
more precisely in Section 3.}
such as $1, -1$, or $2$), corresponding to the fact that 
$|NC(3)| = 5$.

Let now $( \cA , \varphi , \phiprim )$ be an incps as in 
Definition \ref{def:1.1}. Then in addition to the non-crossing 
cumulant functionals $\kappa_n : \cA^n \to \bC$ associated to 
$\varphi$ we will define another family of multilinear 
functionals $( \kappa_n ' : \cA^n \to \bC )_{n \geq 1}$, which
involve both $\varphi$ and $\phiprim$. For every $n \geq 1$, the 
functional $\kappa_n '$ is obtained by taking a 
{\em formal derivative} in the formula for $\kappa_n$, where we 
postulate that the 
derivative of $\varphi$ is $\phiprim$ and we invoke linearity
and the Leibnitz rule for derivatives. For instance for $n=3$
the term $\varphi (a_1 a_2 a_3)$ on the right-hand side of 
(\ref{eqn:1.21}) is derivated into $\phiprim (a_1 a_2 a_3)$,
the term $\varphi (a_1) \varphi (a_2 a_3)$ is derivated into 
$\phiprim (a_1) \varphi (a_2 a_3) +
\varphi (a_1) \phiprim (a_2 a_3)$, etc, yielding the formula 
for $\kappa_3 '$ to be
\begin{equation}    \label{eqn:1.22}
\kappa_3 ' (a_1, a_2, a_3) = 
\begin{array}[t]{l}
\phiprim (a_1 a_2 a_3) - \phiprim (a_1) \varphi (a_2 a_3)  
                       - \varphi (a_1) \phiprim (a_2 a_3)   \\
- \phiprim (a_2) \varphi (a_1 a_3)   
                       - \varphi (a_2) \phiprim (a_1 a_3)   
- \phiprim (a_3) \varphi (a_1 a_2)
                       - \varphi (a_3) \phiprim (a_1 a_2)   \\
+ 2 \phiprim (a_1) \varphi (a_2) \varphi (a_3)
+ 2 \varphi (a_1) \phiprim (a_2) \varphi (a_3)
+ 2 \varphi (a_1) \varphi (a_2) \phiprim (a_3).
\end{array}
\end{equation}  
We will refer to the functionals $\kappa_n '$ as 
{\em infinitesimal non-crossing cumulants} associated
to $( \cA , \varphi , \phiprim )$. The precise formula defining 
them appears in Definition \ref{def:4.2} below.
The passage from the formula for $\kappa_n$ to the one for
$\kappa_n '$ is related to a concept of {\em dual derivation
system} on a space of multilinear functionals on $\cA$, which 
is discussed in Section 7 of the paper.

The role of infinitesimal non-crossing cumulants in the study 
of infinitesimal freeness is described in the next theorem.

\begin{theorem}    \label{thm:1.2}
Let $( \cA , \varphi , \phiprim )$ be an incps and let 
$\cA_1, \ldots , \cA_k$ be unital subalgebras of $\cA$. The 
following statements are equivalent:

\noindent
(1) $\cA_1, \ldots , \cA_k$ are infinitesimally free.

\noindent
(2) For every $n \geq 2$, for every 
$i_1, \ldots , i_n \in \{ 1, \ldots , k \}$ which are not 
all equal to each other, and for every
$a_1 \in \cA_{i_1}, \ldots , a_n \in \cA_{i_n}$, one has that 
$\kappa_n (a_1, \ldots , a_n) = \kappa_n ' (a_1, \ldots , a_n) = 0$.
\end{theorem}

Theorem \ref{thm:1.2} provides an infinitesimal version for the 
basic result of Speicher which describes the usual freeness of 
$\cA_1, \ldots , \cA_k$ in $( \cA , \varphi )$ in terms of the 
cumulants $\kappa_n$ (cf. \cite{NS2006}, Theorem 11.16).

In the remaining part of this subsection we point out some other
interpretations of the formula defining $\kappa_n '$ (all  
corresponding to one or another of the points of view (a), (b), (c) 
listed at the beginning of subsection 1.1). The easy verifications
required by these alternative descriptions of $\kappa_n '$ are 
shown at the beginning of Section 4.

First of all one can consider, as in \cite{BS2009}, the situation
when $\varphi , \phiprim$ in (\ref{eqn:1.11}) are obtained as the 
{\em infinitesimal limit} of a family of functionals
$\{ \varphi_t \mid t \in T \}$. Here $T$ is a subset
of $\bR$ which has $0$ as an accumulation point, every $\varphi_t$
is linear with $\varphi_t (1_{\cA}) = 1$, and we have 
\begin{equation}         \label{eqn:1.23}
\varphi (a) = \lim_{t \to 0} \varphi_t (a) \ \mbox{ and }
\ \phiprim (a) = \lim_{t \to 0} 
\frac{\varphi_t (a) - \varphi (a)}{t}, \ \ 
\forall \, a \in \cA.
\end{equation}    
(Note that such families $\{ \varphi_t \mid t \in T \}$ can 
in fact always be found, e.g. by simply taking 
$\varphi_t = \varphi + t \phiprim$, $t \in ( 0, \infty )$.)
In such a situation, the formal derivative which leads from 
$\kappa_n$ to $\kappa_n '$ turns out to have the same effect as
a ``$\frac{d}{dt}$'' derivative. Consequently, we get the 
alternative formula
\begin{equation}   \label{eqn:1.24}
\kappa_n ' (a_1, \ldots , a_n) =
\Bigl[ \ \frac{d}{dt} \kk^{(t)}_n (a_1, \ldots , a_n) \, \Bigr]
\ \vline \ { }_{ { }_{t=0} },
\end{equation}
where $\kk^{(t)}_n$ denotes the $n$th non-crossing 
cumulant functional of $\varphi_{t}$. 

Second of all, it is possible to take a direct combinatorial approach 
to the functionals $\kappa_n '$, and identify precisely a set of 
non-crossing partitions which indexes the terms in the summation 
defining $\kappa_n ' (a_1, \ldots , a_n )$. This set turns out to be 
\begin{equation}   \label{eqn:1.25}
NCZ^{(B)} (n) := \{ \tau \in \ncb (n) \mid \tau
\mbox{ has a zero-block} \} ,
\end{equation}  
where $\ncb (n)$ is the lattice of non-crossing partitions of 
type B of $\{ 1, \ldots , n \} \cup \{ -1, \ldots , -n \}$
(see subsection 3.1 for a brief review of this). Hence in a 
terminology focused on types of non-crossing partitions, 
one could call the functionals $\kappa_n$ and $\kappa_n '$ 
``non-crossing cumulants of type A and of type B'',
respectively. The idea put forth here is that, in some 
sense, summations over $\nczb (n)$ appear as ``derivatives 
for summation over $NC(n)$''. A more refined formula 
supporting this idea is shown in Proposition \ref{prop:7.5} 
below, in connection to the concept of dual derivation sytem.

In the case $n=3$ that we are using for illustration, the 
10 terms appearing on the right-hand side of (\ref{eqn:1.22}) are 
indexed by the 10 partitions with zero-block in $\ncb (3)$. 
For example, the partitions corresponding to the first three 
terms and the last term from (\ref{eqn:1.22}) are depicted in 
Figure 1. The relation between a partition $\tau$ and the 
corresponding term is easy to follow: the zero-block $Z$ of $\tau$ 
produces the $\phiprim ( \, \cdots \, )$ factor, and every pair 
$V, -V$ of non-zero-blocks of $\tau$ produces a 
$\varphi ( \, \cdots \, )$ factor.

\begin{center}
\scalebox{.48}{\includegraphics{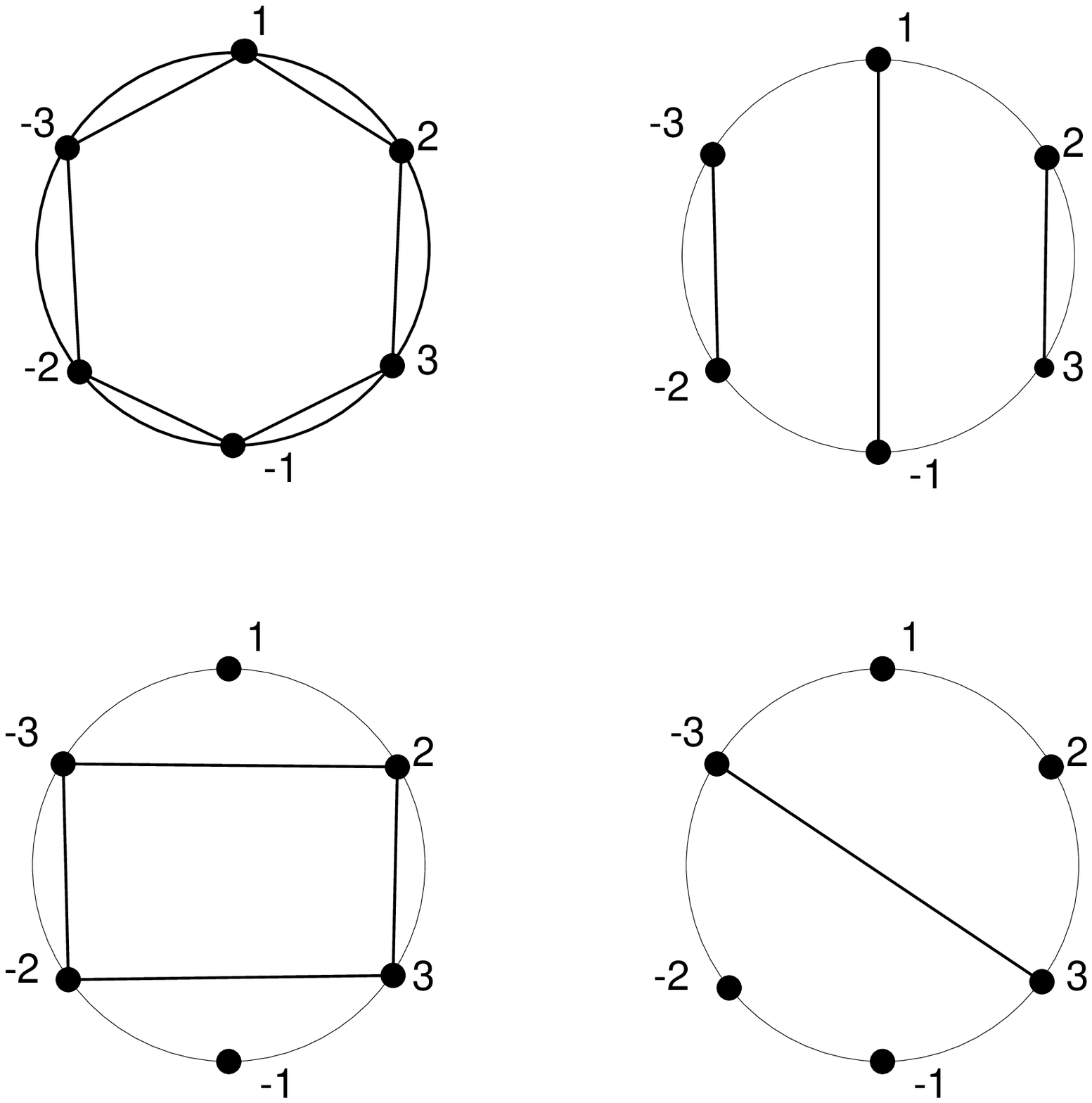}}

\vspace{-3.2cm}

{\bf Figure 1.} Some partitions in $NCZ^{(B)} (3)$.
\end{center}

$\ $

Finally (third of all) one can also give a description of 
$\kappa_n '$ which corresponds to the ``$\bG$-valued'' point of 
view appearing as (b) on the list from subsection 1.1. This goes 
as follows. Let $\phitild = \varphi + \ee \phiprim : \cA \to \bG$ 
be as in (\ref{eqn:1.12}), and consider the family of $\bC$-multilinear 
functionals $( \ktild_n : \cA^n \to \bG )_{n \geq 1}$ defined by 
the same summation formula as for the usual non-crossing cumulant 
functionals $( \kappa_n : \cA^n \to \bC )_{n \geq 1}$, only that 
now we use $\phitild$ instead of $\varphi$ in the summations. So, 
for example, for $n=3$ we have
\begin{equation}    \label{eqn:1.26}
\ktild_3 (a_1, a_2, a_3) = 
\begin{array}[t]{l}
\phitild (a_1 a_2 a_3) - \phitild (a_1) \phitild (a_2 a_3)  
                      - \phitild (a_2) \phitild (a_1 a_3)   \\
- \phitild (a_3) \phitild (a_1 a_2)
+ 2 \phitild (a_1) \phitild (a_2) \phitild (a_3) \in \bG , \ \ 
\ \ \forall \, a_1, a_2, a_3 \in \cA .
\end{array}
\end{equation}  
It then turns out that the functional $\kappa_n '$ can be obtained
by reading the $\ee$-component of $\ktild_n$.

We take the opportunity to introduce here a piece of terminology 
from the literature on Grassman algebras (see e.g. \cite{DW1992}, 
pp. 1-2): the complex numbers $\alpha , \beta$ which give the two 
components of a Grassman number 
$\gamma = \alpha + \ee \beta \in \bG$ will be called the {\em body} 
and respectively the {\em soul} of $\gamma$; it will come in handy 
throughout the paper to denote them 
\footnote{ Besides being amusing, ``Bo'' and ``So'' give a faithful 
analogue for the common notations ``Re'' and ``Im'' used 
when one introduces $\bC$ as a 2-dimensional algebra over $\bR$. }
as
\begin{equation}   \label{eqn:1.27}
\alpha = \bo ( \gamma ), \ \ \beta  = \so ( \gamma ). 
\end{equation}
This notation will also be used in connection to a $\bG$-valued 
function $f$ defined on some set $\cS$ -- we define functions 
$\bo \, f$ and $\so \, f$ from $\cS$ to $\bC$ by
\begin{equation}   \label{eqn:1.28}
( \bo \, f ) (x) = \bo \, (f(x)), 
\ \ ( \so \, f ) (x) = \so \, (f(x)), \ \ \forall \, x \in \cS .
\end{equation}
Returning then to the functionals $\ktild_n : \cA^n \to \bG$ from
the preceding paragraph, their connection to the $\kappa_n '$ (and
also to the $\kappa_n$) can be recorded as
\begin{equation}    \label{eqn:1.29}
\bo \ \ktild_n  = \kappa_n, \ \ 
\so \ \ktild_n  = \kappa_n ' , \ \ \forall \, n \geq 1.
\end{equation}
Due to (\ref{eqn:1.29}), $\ktild_n$ can be used as a simplifying 
tool in calculations with $\kappa_n '$ (in the sense that it may 
be easier to run the corresponding calculation with $\ktild_n$, in 
$\bG$, and only pick soul parts at the end of the calculation). In 
particular, this will be useful when proving Theorem \ref{thm:1.2},
since the condition $\kappa_n (a_1, \ldots , a_n)$ = 
$\kappa_n ' (a_1, \ldots , a_n) = 0$ from 
Theorem \ref{thm:1.2}(2) amounts precisely to
$\ktild_n (a_1, \ldots , a_n) = 0$.

$\ $

\begin{center}
{\bf 1.3 Using derivations to find ``soul companions'' for a 
given \boldmath{$\varphi$} }
\end{center}

When studying infinitesimal freeness it may be of interest 
to consider the situation where we have fixed a noncommutative 
probability space $( \cA , \varphi )$ and a family 
$\cA_1, \ldots , \cA_k$ of unital subalgebras of $\cA$ which 
are free in $( \cA , \varphi )$. In this situation we can ask: 
how do we find interesting examples of functionals 
$\phiprim : \cA \to \bC$ with $\phiprim ( 1_{\cA} ) = 0$ and 
such that $\cA_1, \ldots , \cA_k$ become infinitesimally free in 
$( \cA , \varphi , \phiprim )$? A nice name for such 
functionals $\phiprim$ is suggested by the $\bG$-valued point 
of view described in subsection 1.1: since 
$\varphi$ and $\phiprim$ are the body part and respectively the 
soul part of the consolidated functional $\phitild : \cA \to \bG$, 
one may say that we are looking for a suitable 
{\em soul companion} $\phiprim$ for the given ``body functional'' 
$\varphi$ (and in reference to the given subalgebras 
$\cA_1, \ldots , \cA_k$).

Let us note that the remark made at the end of subsection 1.1 can
be interpreted as a statement about soul companions. Indeed, this 
remark says that if $( \cA , \varphi )$ is the free product of 
$( \cA_1, \varphi_1 ), \ldots , ( \cA_k , \varphi_k )$, then a
$\phiprim$ from the desired set of soul companions is parametrized 
precisely by a family of linear functionals 
$\phiprim_i : \cA_i \to \bC$ such that $\phiprim_i ( 1_{\cA} ) = 0$, 
$1 \leq i \leq k$.

The point we follow here, with inspiration from \cite{BS2009}, is 
that some interesting recipes to construct ``soul companions'' 
for a given $\varphi : \cA \to \bC$ arise from ideas pertaining to 
differentiability. This is intimately related to the fact that 
$\kappa_n '$ is a formal derivative for $\kappa_n$, hence to 
equations of the form 
\[
d_n ( \kappa_n ) = \kappa_n ', \ \ \forall \, n \geq 1,
\]
where $(d_n)_{n \geq 1}$ is a dual derivation system on $\cA$.
Indeed, suppose we are given a derivation $D : \cA \to \cA$; then one 
has a natural dual derivation system associated to it, which acts by
\begin{equation}   \label{eqn:1.31}
( d_n f ) (a_1, \ldots , a_n) = \sum_{m=1}^n 
f \bigl( a_1, \ldots , a_{m-1}, D(a_m), a_{m+1}, \ldots , a_n \bigr),
\end{equation}
for $f : \cA^n \to \bC$ multilinear and  $a_1, \ldots , a_n \in \cA$.
By using the $d_n$ from (\ref{eqn:1.31}), we obtain the following 
theorem.

\begin{theorem}    \label{thm:1.3}
Let $( \cA ,\varphi , \phiprim )$ be an incps, and let 
$\kk_n$ and $\kk_n '$ be the non-crossing cumulant functionals
associated to it. Suppose $D : \cA \to \cA$ is a derivation 
with the property that $\phiprim =\varphi \circ D$. Then for 
every $n \geq 1$ and $a_1, \ldots , a_n \in \cA$ one has
\begin{equation}  \label{eqn:1.32}  
\kappa_n '  (a_1, \ldots , a_n) = \sum_{m=1}^n 
\kappa_n (a_1, \ldots , a_{m-1}, D(a_m), a_{m+1}, \ldots ,a_n ).
\end{equation}
\end{theorem}

Moreover, when combined with Theorem \ref{thm:1.2}, the formula 
for infinitesimal cumulants obtained in (\ref{eqn:1.32}) has the 
following immediate consequence.

\begin{corollary}    \label{cor:1.4}
Let $( \cA , \varphi )$ be a noncommutative probability space, and 
let $\cA_1, \ldots , \cA_k$ be unital subalgebras of $\cA$ which 
are free in $( \cA , \varphi )$. Suppose we found a derivation 
$D: \cA \to \cA$ such that $D( \cA_i ) \subseteq \cA_i$ for every 
$1 \leq i \leq k$.  Then $\cA_1, \ldots , \cA_k$ are infinitesimally 
free in $( \cA , \varphi , \phiprim )$, where 
$\phiprim = \varphi \circ D$.
\end{corollary}

For comparison, let us also look at the parallel statement arising
in connection to infinitesimal limits. This is essentially the same 
as Remark 15 from \cite{BS2009}, and goes as follows.

\begin{proposition}    \label{prop:1.5}
Let $( \cA , \varphi )$ be a noncommutative probability space, and 
let $\cA_1, \ldots , \cA_k$ be unital subalgebras of $\cA$ which 
are free in $( \cA , \varphi )$. Suppose we found a family of 
linear functionals $( \, \varphi_t : \cA \to \bC \, )_{t \in T}$
with $\varphi_t ( 1_{\cA} ) = 1$ for every $t \in T$ and
such that:

(i) $\cA_1, \ldots , \cA_k$ are free in $( \cA , \varphi_t )$ 
for every $t \in T$.

(ii) $\lim_{t \to 0} \varphi_t (a) = \varphi (a)$, for every
$a \in \cA$.

(iii) The limit $\phiprim (a) := \lim_{t \to 0}
( \varphi_t (a) - \varphi (a) )/t$ exists, for every $a \in \cA$.

\noindent
Then $\cA_1, \ldots , \cA_k$ are infinitesimally 
free in $( \cA , \varphi , \phiprim )$, where 
$\phiprim : \cA \to \bC$ is defined by condition (iii).
\end{proposition}

A natural example accompanying Proposition \ref{prop:1.5} comes 
in connection to $\boxplus$-convolution powers of joint 
distributions of $k$-tuples (cf. Example \ref{ex:8.9} below).
In Section 8 we also discuss a couple of natural situations 
when Corollary \ref{cor:1.4} applies (cf. Example \ref{ex:8.7}).

$\ $

\begin{center}
{\bf 1.4 Outline of the rest of the paper}
\end{center}

Besides the introduction, the paper has seven other sections. 
In Section 2 we collect some basic properties of infinitesimal 
freeness, and we discuss the relations between Definition 
\ref{def:1.1} and the frameworks of \cite{BS2009}, \cite{BGN2003}.
Section 3 is a review of background concerning non-crossing 
partitions and non-crossing cumulants. In Section 4 we introduce
the non-crossing infinitesimal cumulants, we verify the 
equivalence between their various alternative descriptions, and 
we prove Theorem \ref{thm:1.2}. 

Sections 5 and 6 address the topic of alternating products of 
infinitesimally free random variables. Section 5 uses this topic 
to illustrate a ``generic'' method to obtain infinitesimal 
analogues for known results in usual free probability: one 
replaces $\bC$ by $\bG$ in the proof of the original result, then 
one takes the soul part in the $\bG$-valued statement that comes 
out. By using this method we obtain the infinitesimal versions of 
two important facts related to alternating products that were 
originally found in \cite{NS1996} -- one of them is about 
compressions by free projections, the other concerns a method 
of constructing free families of free Poisson elements. In 
Section 6 we remember that the concept of incps has its origins 
in the considerations ``of type B'' from \cite{BGN2003}, and we 
look at how the essence of these considerations persists in the 
framework of the present paper. The main point of the section 
is that, when taking the soul part of the $\bG$-valued formulas
for alternating products of infinitesimally free random variables,
one does indeed obtain nice analogues of type B (with summations 
over $\ncb (n)$) for the type A formulas. In particular, this 
offers another explanation for why the infinitesimal cumulant 
functional $\kappa_n '$ can be described by using a summation 
formula over $\nczb (n)$.

In Section 7 we return to the point of view of treating 
$\kappa_n '$ as a derivative of the usual non-crossing
cumulant functional $\kappa_n$, and we discuss the related 
concept of dual derivation system on a unital algebra $\cA$.
Finally, Section 8 elaborates on the discussion about soul
companions from the above subsection 1.3. In particular, we 
show how the dual derivation system provided by a derivation 
$D : \cA \to \cA$ leads to the setting for infinitesimal 
freeness from Corollary \ref{cor:1.4}. Section 8 (and the paper) 
concludes with a couple of examples related to the settings of 
Corollary \ref{cor:1.4} and of Proposition \ref{prop:1.5}. 

$\ $

$\ $

\begin{center}
{\bf\large 2. Basic properties of infinitesimal freeness}
\end{center}
\setcounter{section}{2}
\setcounter{equation}{0}
\setcounter{theorem}{0}

In this section we collect some basic properties of infinitesimal
freeness, and we discuss the relations between Definition 
\ref{def:1.1} and the frameworks from \cite{BS2009}, 
\cite{BGN2003}.

\begin{definition}     \label{def:2.1}
Here are some standard variations of Definition \ref{def:1.1}.

$1^o$ The concept of infinitesimal freeness carries 
over to $*$-algebras. More precisely, we will use the name
{\em $*$-incps} for an incps $( \cA , \varphi , \phiprim )$ where 
$\cA$ is a unital $*$-algebra and where 
 
(i) $\varphi$ is positive definite, that is,
$\varphi (a^* a) \geq 0, \ \forall \, a \in \cA$;
 
(ii) $\phiprim$ is selfadjoint, that is,
$\phiprim (a^*) = \overline{\phiprim (a)}, \ \forall \, a \in \cA$.

$2^o$ Another standard variation of the definitions is that 
infinitesimal freeness can be considered for arbitrary subsets 
of $\cA$ (which don't have to be subalgebras). So if 
$( \cA , \varphi , \phiprim )$ is an incps (respectively a 
$*$-incps) and if $\cX_1, \ldots , \cX_k$ are subsets of $\cA$, 
then we will say that $\cX_1, \ldots , \cX_k$ are 
{\em infinitesimally free} (respectively 
{\em infinitesimally $*$-free}) when the unital subalgebras 
(respectively $*$-subalgebras) generated by 
$\cX_1, \ldots , \cX_k$ are so.
\end{definition}

\begin{remark}    \label{rem:2.2}
Let $( \cA , \varphi )$ be a noncommutative probability space and 
let $\cA_1 , \ldots , \cA_k$ be unital subalgebras of $\cA$ which 
are free in $( \cA , \varphi )$. It is very easy to see 
(cf. Remark 2.5.2 in \cite{VDN1992} or Examples 5.15 in 
\cite{NS2006}) that the way how $\varphi$ acts on 
$\mbox{Alg} ( \cA_1 \cup \cdots \cup \cA_k )$ 
can be reconstructed from the restrictions $\varphi \mid \cA_i$, 
$1 \leq i \leq k$. The simplest illustration for how this works
is provided by the formula 
\begin{equation}  \label{eqn:2.21}
\varphi (ab) = \varphi (a) \varphi (b),
\ \ \forall \, a \in \cA_{i_1}, \, b \in \cA_{i_2}, 
\mbox{ with $i_1 \neq i_2$,}
\end{equation}
which is obtained by expanding the product and then collecting 
terms in the equation
$\varphi \Bigl( \, (a - \varphi (a) 1_{\cA}) \cdot 
(b - \varphi (b) 1_{\cA}) \, \Bigr) = 0$.

A similar phenomenon turns out to take place
when dealing with infinitesimal freeness: the way how
$\phiprim$ acts on $\mbox{Alg} ( \cA_1 \cup \cdots \cup \cA_k )$ 
can be reconstructed from the restrictions of $\varphi$ and of 
$\phiprim$ to $\cA_i$, $1 \leq i \leq k$. For example, the 
counterpart of Equation (\ref{eqn:2.21}) says that
\begin{equation}  \label{eqn:2.22}
\phiprim (ab) = \phiprim (a) \varphi (b) + \varphi (a) \phiprim (b),
\ \ \forall \, a \in \cA_{i_1}, \, b \in \cA_{i_2}, 
\mbox{ where $i_1 \neq i_2$.}
\end{equation}
This is obtained by expanding the product and then collecting 
terms in the equation
$\phiprim \Bigl( \, (a - \varphi (a) 1_{\cA}) \cdot 
(b - \varphi (b) 1_{\cA}) \, \Bigr) = 0$ (which is a particular
case of Equation (\ref{eqn:1.15})), and by taking into account 
that $\phiprim ( 1_{\cA} ) = 0$.

We leave it as an easy exercise to the reader to verify that the 
similar calculation for an alternating product of 3 factors (which
makes a more involved use of Equation (\ref{eqn:1.15})) leads 
to the formula
\begin{equation}  \label{eqn:2.23}
\phiprim (a_1 b a_2) = \phiprim (a_1 a_2) \varphi (b) + 
\varphi (a_1 a_2) \phiprim (b),
\ \ \mbox{ for } a_1, a_2 \in \cA_{i_1}, \, b \in \cA_{i_2}, 
\mbox{ with $i_1 \neq i_2$.}
\end{equation}
\end{remark}

\begin{remark}     \label{rem:2.3}
{\em (Traciality.)} Another well-known fact in usual free probability 
is that if the unital subalgebras 
$\cA_1, \ldots , \cA_k \subseteq \cA$ are free in $( \cA , \varphi )$ 
and if $\varphi \mid \cA_i$ is a trace for every $1 \leq i \leq k$, 
then $\varphi$ is a trace on 
$\mbox{Alg} ( \cA_1 \cup \cdots \cup \cA_k )$. This too extends to 
the infinitesimal framework: if $\cA_1, \ldots , \cA_k$ are 
infinitesimally free in $( \cA , \varphi , \phiprim )$ and if 
$\varphi \mid \cA_i, \phiprim \mid \cA_i$ are traces for 
every $1 \leq i \leq k$, then $\varphi$ and $\phiprim$ are 
traces on Alg$( \cA_1 \cup \cdots \cup \cA_k )$. Rather than 
writing an ad-hoc proof of this fact based directly on Definition 
\ref{def:1.1}, we find it more instructive to do this by using 
cumulants -- see Proposition \ref{prop:4.10} below.

We next move to describing the free product of infinitesimal 
noncommutative probability spaces announced at the end of 
Section 1.1. 
\end{remark}

\begin{proposition}     \label{prop:2.4}
Let $( \cA_1 , \varphi_1 ), \ldots , ( \cA_k , \varphi_k )$ be
noncommutative probability spaces, and consider the free 
product $( \cA , \varphi ) = 
( \cA_1 , \varphi_1 ) * \cdots * ( \cA_k , \varphi_k )$
(as described e.g. in Lecture 6 of \cite{NS2006}). Suppose that 
for every $1 \leq i \leq k$ we are given a linear functional 
$\phiprim_i : \cA_i \to \bC$ such that $\phiprim_i ( 1_{\cA} ) = 0$. 
Then there exists a unique linear functional $\phiprim : \cA \to \bC$ 
such that $\phiprim \mid \cA_i = \phiprim_i$, $1 \leq i \leq k$, 
and such that $\cA_1, \ldots , \cA_k$ are infinitesimally free in 
$( \cA , \varphi , \phiprim )$.
\end{proposition}

\begin{proof} We start by reviewing a few basic facts and notations 
related to $( \cA , \varphi )$. Each of $\cA_1 , \ldots , \cA_k$ 
is identified as a unital subalgebra of $\cA$, such that 
$\varphi \mid \cA_i = \varphi_i$. For $1 \leq i \leq k$ we denote
$\cAo_i = \{ a \in \cA_i \mid \varphi (a) = 0 \}$,
and for every $n \geq 1$ and $1 \leq i_1, \ldots , i_n \leq k$ 
such that $i_1 \neq i_2, \ldots , i_{n-1} \neq i_n$ we put
\begin{equation}    \label{eqn:2.41}
\cW_{i_1, \ldots , i_n} := \mbox{span} \bigl\{ a_1 \cdots a_n 
\mid a_1 \in \cAo_{i_1}, \ldots , a_n \in \cAo_{i_n} \bigr\} .
\end{equation} 
It is known that $\cW_{i_1, \ldots , i_n}$ is canonically 
isomorphic to the tensor product 
$\cAo_{i_1} \otimes \cdots \otimes \cAo_{i_n}$, via the 
identification 
$a_1 \cdots a_n \simeq a_1 \otimes \cdots \otimes a_n$,
for $a_1 \in \cAo_{i_1}, \ldots , a_n \in \cAo_{i_n}$.
Moreover it is known that the spaces $\cW_{i_1, \ldots , i_n}$ 
defined in (\ref{eqn:2.41}) realize a direct sum decomposition of 
the kernel of $\varphi$. (See \cite{NS2006}, pp. 81-84.)

Due to the direct sum decomposition mentioned above, we may define
the required functional $\phiprim$ by separately prescribing its 
behaviour at $1_{\cA}$ and on each of the subspaces 
$\cW_{i_1, \ldots , i_n}$. We put $\phiprim (1_{\cA}) := 0$. We 
also prescribe $\phiprim$ to be $0$ on $\cW_{i_1, \ldots , i_n}$
whenever $n$ is even, and whenever $n$ is odd but it is not true
that $i_m = i_{n+1-m}$ for all $1 \leq m \leq (n-1)/2$.
Suppose next that $n = 2m-1$, odd, and that the indices 
$i_1, \ldots , i_n$ are such that  
$i_1 = i_{2m-1}, \, i_2 = i_{2m-2}, \ldots , i_{m-1} = i_{m+1}$.
By using the identification 
$\cW_{i_1, \ldots , i_n} \simeq
\cAo_{i_1} \otimes \cdots \otimes \cAo_{i_n}$
it is immediate that we can 
define a linear map on $\cW_{i_1, \ldots , i_n}$ by the 
requirement that 
\[
a_1 \cdots a_{2m-1}  \mapsto
\varphi_{i_1} (a_1 a_{2m-1}) \varphi_{i_2} (a_2 a_{2m-2}) \cdots  
\varphi_{i_{m-1}} (a_{m-1} a_{m+1}) \cdot  \phiprim_{i_m} (a_m),
\]
for every $a_1 \in \cAo_{i_1}, \ldots , a_n \in \cAo_{i_n}$;
we take this as the prescription for how $\phiprim$ is to act on
$\cW_{i_1, \ldots , i_n}$. 

Directly from Definition \ref{def:1.1} it is immediate that, with 
$\phiprim : \cA \to \bC$ defined as in the preceding paragraph, 
$\cA_1 , \ldots , \cA_k$ are infinitesimally free in 
$( \cA , \varphi , \phiprim )$. The uniqueness of $\phiprim$ with 
this property is also immediate.
\end{proof}

\begin{definition}    \label{def:2.5}
Let $( \cA_1 , \varphi_1, \phiprim_1), \ldots ,
( \cA_k , \varphi_k, \phiprim_k)$ be infinitesimal noncommutative
probability spaces. We define their {\em free product} to be 
$( \cA , \varphi , \phiprim )$ where $( \cA , \varphi )$ = 
$( \cA_1 , \varphi_1 ) * \cdots * ( \cA_k , \varphi_k )$
and where $\phiprim : \cA \to \bC$ is the functional provided by 
Proposition \ref{prop:2.4}.
\end{definition}

\begin{remark}    \label{rem:2.6}
In the context of Proposition \ref{prop:2.4}, suppose that 
$( \cA_1 , \varphi_1 ), \ldots , ( \cA_k , \varphi_k )$ are 
$*$-probability spaces. Then so is the free product 
$( \cA , \varphi )$ (see \cite{NS2006}, Theorem 6.13). If
moreover each of the functionals $\phiprim_i : \cA_i \to \bC$ 
given in Proposition \ref{prop:2.4} is selfadjoint, then it 
is easily checked that the resulting functional 
$\phiprim : \cA \to \bC$ is selfadjoint too. Hence if in 
Definition \ref{def:2.5} each of 
$( \cA_i , \varphi_i , \phiprim_i )$ is a $*$-incps, then the 
free product $( \cA , \varphi , \phiprim )$ is a $*$-incps as
well.
\end{remark}

\begin{example}    \label{ex:2.7}
For an illustration of the above, we look at a simple 
example where the spaces $\cW_{i_1, \ldots , i_n}$ are 
all 1-dimensional. Consider the $k$-fold free product group 
$\bZ_2 * \cdots * \bZ_2$ and let $\varphi$ be the canonical trace 
on the group algebra $\cA := \bC [ \bZ_2 * \cdots * \bZ_2 ]$. So 
$\cA$ is a unital $*$-algebra freely generated by $k$ unitaries 
$u_1, \ldots ,u_k$ of order 2, and has a linear basis $\cB$ given by
\begin{equation}    \label{eqn:2.71}
\cB = \{ 1_{\cA} \} \cup \Bigl\{ u_{i_1} \cdots u_{i_n} \
\begin{array}{ll}
\vline & n \geq 1, \ 1 \leq i_1, \ldots , i_n \leq k, \\
\vline & \mbox{ with } i_1 \neq i_2, \ldots , i_{n-1} \neq i_n  
\end{array}  \Bigr\} .
\end{equation}
The linear functional $\varphi : \cA \to \bC$ acts on the basis 
$\cB$ by 
\[
\varphi ( 1_{\cA} ) = 1, \ \ \mbox{ and } \ \
\varphi ( b) = 0, \ \ \forall \, b \in \cB \setminus \{ 1_{\cA} \}.
\]
It is easy to verify (see e.g. Lecture 6 in \cite{NS2006}) that 
we have $( \cA , \varphi )$ = 
$( \cA_1 , \varphi_1 ) * \cdots * ( \cA_k , \varphi_k )$, 
where for $1 \leq i \leq k$ we denote 
$\cA_i = \mbox{span} \{ 1_{\cA}, u_i \}$
(2-dimensional $*$-subalgebra of $\cA$), and where 
$\varphi_i := \varphi \mid \cA_i$. The 
direct sum decomposition of $\cA$ with respect to this free product
structure simply has 
\[
\cW_{i_1, \ldots , i_n} = \mbox{$1$-dimensional space spanned by }
u_{i_1} \cdots u_{i_n},
\]
for every $n \geq 1$ and every alternating sequence 
$i_1, \ldots , i_n$ as described in (\ref{eqn:2.71}).

Now let $\phiprim_i : \cA_i \to \bC$ be linear functionals such 
that $\phiprim_i ( 1_{\cA} ) = 0$, $1 \leq i \leq k$. Clearly, these 
functionals are determined by the values 
\[
\phiprim_1 ( u_1 ) =: \alpha_1 ' , \ldots ,
\phiprim_k ( u_k ) =: \alpha_k ' . 
\]
The free product extension $\phiprim : \cA \to \bC$ then acts by
\begin{equation}   \label{eqn:2.72}
\phiprim (u_{i_1} \cdots u_{i_n}) = 
\left\{  \begin{array}{ll}
\alpha_{i_m} ', & \mbox{ if $n$ is odd, $n = 2m-1$, and 
             $i_1 = i_{2m-1}, \ldots , i_{m-1} = i_{m+1}$}  \\
0,              & \mbox{ otherwise.}
\end{array}  \right.
\end{equation}
Note that formula (\ref{eqn:2.72}) looks particularly nice in the 
case when $k=2$ -- indeed, in this case the requirement that 
$i_1 = i_{2m-1}, \ldots , i_{m-1} = i_{m+1}$ is automatically 
satisfied whenever $n = 2m-1$ and $i_1, \ldots , i_n$ are as in
(\ref{eqn:2.71}).
\end{example}

\begin{remark}      \label{rem:2.8}
{\em (Relation to \cite{BS2009}).} 
Definition 13 of \cite{BS2009} introduces a concept of 
infinitesimal freeness for unital subalgebras 
$\cA_1, \cA_2 \subseteq \cA$ in an incps $( \cA , \mu , \mu ' )$.
As explained there (immediately following to Definition 13), this
amounts to two requirements: that $\cA_1, \cA_2$ are free 
in $( \cA , \mu )$, and that they satisfy the following additional 
condition:
\begin{equation}  \label{eqn:2.81}
\mu' \Bigl( \, \bigl( p_1-\mu (p_1) 1_{\cA} \bigr) \cdots 
\bigl( p_n-\mu (p_n) 1_{\cA} \bigr) \, \Bigr) = 
\end{equation}
\[
\sum_{m=1}^n \mu \Bigl( \, (p_1- \mu (p_1) 1_{\cA}) \cdots 
\mu '(p_m) \cdots (p_n-\mu (p_n) 1_{\cA}) \, \Bigr)
\]
for $p_1\in \cA_{i_1}, \ldots ,p_n\in \cA_{i_n}$, 
where $i_1 \neq i_2, \ldots , i_{n-1} \neq i_n$. By denoting 
$p_m - \mu (p_m) 1_{\cA} =: q_m$ and by taking into account 
that $\mu ' (q_m) = \mu ' (p_m)$, $1 \leq m \leq n$,
one sees that condition (\ref{eqn:2.81}) is equivalent to its 
particular case requesting that 
\begin{equation}  \label{eqn:2.82}
\mu' (q_1 \cdots q_n) = \sum_{m=1}^n 
\mu (q_1 \cdots q_{m-1} q_{m+1} \cdots q_n) \cdot \mu '(q_m)
\end{equation}
for $q_1\in \cA_{i_1}, \ldots ,q_n\in \cA_{i_n}$, 
where $i_1 \neq i_2, \ldots , i_{n-1} \neq i_n$ and where 
$\mu (q_1) = \cdots = \mu (q_n) = 0$. 

But now, let $\cA_1 , \cA_2$ be unital subalgebras of $\cA$ 
which are free in $( \cA , \mu )$. A standard 
calculation from usual free probability (see e.g. Lemma 5.18 on 
page 73 of \cite{NS2006}) says that, with $q_1, \ldots , q_n$ as 
in (\ref{eqn:2.82}), one has 
$\mu (q_1 \cdots q_{m-1} q_{m+1} \cdots q_n) = 0$ unless it 
is true that $m-1 = n-m$ and that 
$i_{m-1} = i_{m+1}, i_{m-2} = i_{m+2}, \ldots , i_1 = i_n$; 
moreover, if the latter conditions are satisfied, then 
\[
\mu (q_1 \cdots q_{m-1} q_{m+1} \cdots q_n) = 
\mu (q_{m-1} q_{m+1}) \,
\mu (q_{m-2} q_{m+2}) \cdots \mu (q_1 q_n).
\]
This clearly implies that the sum on the right-hand side of 
(\ref{eqn:2.82}) has at most one term which is different from 0;
and moreover, when such a term exists, it is exactly as described 
in Equation (\ref{eqn:1.15}) of Definition \ref{def:1.1}.

Hence, modulo an immediate reformulation, the concept of 
infinitesimal freeness from \cite{BS2009} is the same as the one
used in this paper (which justifies the fact that we are calling 
it by the same name).
\end{remark}

\begin{remark}        \label{rem:2.9}
{\em (Relation to \cite{BGN2003}).}
A {\em noncommutative probability space of type B} is defined in 
\cite{BGN2003} as a system $( \cA , \varphi , \cV , f, \Phi )$,
where $( \cA , \varphi )$ is a noncommutative probability space,
$\cV$ is a complex vector space, $f : \cV \to \bC$ is 
a linear functional, and $\Phi : \cA \times \cV \times \cA \to \cV$
is a two-sided action. We will write for short $a \xi b$ and 
respectively $a \xi$, $\xi b$ instead of 
$\Phi (a, \xi , b)$ and respectively 
$\Phi (a, \xi , 1_{\cA} )$, $\Phi ( 1_{\cA}, \xi , b )$, for 
$a, b \in \cA$ and $\xi \in \cV$. Let $\cA_1, \ldots , \cA_k$ be 
unital subalgebras of $\cA$ and let $\cV_1, \ldots ,  \cV_k$ be 
linear subspaces of $\cV$, such that $\cV_i$ is closed under the 
two-sided action of $\cA_i$, $1 \leq i \leq k$. Definition 7.2 of 
\cite{BGN2003} introduces a concept of what it means for 
$( \cA_1, \cV_1 ), \ldots , ( \cA_k, \cV_k )$ to be free in 
$( \cA , \varphi , \cV , f, \Phi )$. This amounts to two 
requirements: that $\cA_1, \ldots , \cA_k$ are free in 
$( \cA , \varphi )$, and that the following additional condition 
is satisfied:
\begin{equation}  \label{eqn:2.91}
f( a_m \dots a_1 \xi b_1 \dots b_n)=
\left\{   \begin{array}{l}
\varphi (a_1 b_1) \cdots \varphi (a_n b_n) f (\xi),         \\
\mbox{$\ \ $ if $m=n$ and $i_1 = j_1, \ldots , i_n = j_n$}  \\
0, \mbox{ otherwise,}
\end{array}  \right.
\end{equation}
holding for $m, n \geq 0$ and
$a_1 \in \cA_{i_1}, \ldots ,a_m \in \cA_{i_m},
b_1 \in \cA_{j_1}, \ldots ,b_n\in \cA_{j_n}, \xi \in \cV_h$,
where any two consecutive indices among 
$i_m,  \ldots, i_1, h, j_1, \ldots , j_n$ are different from 
each other, and where  $\varphi (a_m) = \cdots = \varphi (a_1)=0$
= $\varphi (b_1)= \cdots =\varphi (b_n)$.

Now, to $( \cA , \varphi , \cV , f, \Phi )$ as above one 
associates a {\em link-algebra}, which is simply the direct 
product $\cM = \cA \times \cV$ endowed with the natural 
structure of complex vector space and with multiplication
\begin{equation}  \label{eqn:2.92}
(a, \xi ) \cdot (b, \eta ) = (ab, a \eta + \xi b ), \ \ 
\forall \, a,b \in \cA, \, \xi , \eta \in \cV .
\end{equation}
If we define $\psi , \psiprim : \cM \to \bC$ by
\begin{equation}  \label{eqn:2.93}
\psi   ( \, (a, \xi ) \, ) := \varphi (a) , \ \ 
\psiprim   ( \, (a, \xi ) \, ) := f (\xi) , \ \ 
\forall \, (a , \xi ) \in \cM ,
\end{equation}
then $( \cM , \psi , \psiprim )$ becomes an incps.
Let again $\cA_1, \ldots , \cA_k$ be unital subalgebras of $\cA$ 
and $\cV_1, \ldots ,  \cV_k$ be linear subspaces of $\cV$ such 
that $\cV_i$ is closed under the two-sided action of $\cA_i$, 
$1 \leq i \leq k$. Then $\cM_1 := \cA_1 \times \cV_1, \ldots ,
\cM_k := \cA_k \times \cV_k$ are unital subalgebras of the 
link-algebra $\cM$, and we claim that 
\begin{equation}  \label{eqn:2.94}
\left(  \begin{array}{c}
( \cA_1, \cV_1 ), \ldots , ( \cA_k, \cV_k )                \\
\mbox{ are free in $( \cA , \varphi , \cV , f, \Phi )$,}   \\
\mbox{ in the sense of \cite{BGN2003} }
\end{array} \right) \ \Leftrightarrow \ 
\left(  \begin{array}{c}
\cM_1 , \ldots , \cM_k \mbox{ are free}         \\
\mbox{ in $( \cM , \psi , \psiprim )$, in the}  \\
\mbox{ sense of Definition \ref{def:1.1} }
\end{array} \right) .
\end{equation}  
In order to prove the implication ``$\Leftarrow$'' in 
(\ref{eqn:2.94}), we only have to write
\[
f( a_m \dots a_1 \xi b_1 \dots b_n) =
\psiprim \bigl( \, 
( a_m , 0_\cV ) \cdots ( a_1 , 0_\cV ) \cdot ( 0_\cA , \xi ) 
\cdot ( b_1 , 0_\cV ) \cdots ( b_n , 0_\cV ) \, \bigr) 
\] 
and then invoke Equation (\ref{eqn:1.15}). For the implication 
``$\Rightarrow$'', consider some elements $(a_1 , \xi_1 ) 
\in \cM_{i_1}, \ldots , (a_n , \xi_n ) \in \cM_{i_n}$ where 
$i_1 \neq i_2, \ldots , i_{n-1} \neq i_n$ and where 
$\psi ( \, (a_1 , \xi_1 ) \, ) = \cdots =
\psi ( \, (a_n , \xi_n ) \, ) = 0$ (which just means that 
$\varphi (a_1) = \cdots = \varphi (a_n) = 0$). By using how the 
multiplication on $\cM$ and how $\psiprim$ are defined, we see that 
\begin{equation}  \label{eqn:2.95}
\psiprim \bigl( \, (a_1 , \xi_1 ) \cdots (a_n, \xi_n ) \, \bigr)
= \sum_{m=1}^n f( a_1 \cdots a_{m-1} \xi_m a_{m+1} \cdots a_n).
\end{equation}
But because of (\ref{eqn:2.91}), at most one term in the sum 
on the right-hand side of (\ref{eqn:2.95}) can be different from
$0$; moreover such a term can only occur for $m = (n+1)/2$, if 
($n$ is odd and) $i_1 = i_{2m-1}, \ldots , i_{m-1} = i_{m+1}$. 
Finally, if the latter equalities of indices are satisfied, then 
the unique term left in the sum from (\ref{eqn:2.95}) is
$\varphi (a_1 a_{2m-1}) \cdots \varphi (a_{m-1} a_{m+1}) 
f( \xi_m )$, and the conditions defining the infinitesimal freeness 
of $\cM_1 , \ldots , \cM_k$ in $( \cM , \psi , \psiprim )$ follow.

Hence, by focusing on the link-algebra, one can incorporate the 
freeness of type B from \cite{BGN2003} into the framework of 
this paper.
\end{remark}

$\ $

$\ $

\begin{center}
{\bf\large 3. Background on non-crossing partitions and 
non-crossing cumulants}
\end{center}
\setcounter{section}{3}
\setcounter{equation}{0}
\setcounter{theorem}{0}

\begin{center}
{\bf 3.1 Non-crossing partitions}
\end{center}

\begin{notation}     \label{def:3.1}
We will use the standard conventions of notation concerning 
non-crossing partitions (as they appear for instance in Lecture 9 
of \cite{NS2006}). So for a positive integer $n$ we denote by 
$NC(n)$ the set of all non-crossing partitions of 
$\{ 1, \ldots , n \}$. We vill use the abbreviation ``$V \in \pi$''
for ``$V$ is a block of $\pi$'', and the number of blocks of 
$\pi \in NC(n)$ will be denoted as $| \pi |$. On $NC(n)$ we will 
consider the partial order given by {\em reverse refinement}; that
is, for $\pi , \rho \in NC(n)$ we write ``$\pi \leq \rho$'' to 
mean that every block of $\rho$ is a union of blocks of $\pi$.
The minimal and maximal element of $( NC(n), \leq )$ are denoted
by $0_n$ (the partition of $\{ 1, \ldots , n \}$ into $n$ blocks
of $1$ element each) and respectively $1_n$ (the partition of
of $\{ 1, \ldots , n \}$ into 1 block of $n$ elements). It is 
easy to see that $( NC(n), \leq )$ is a lattice, i.e.  that 
every $\pi , \rho \in NC(n)$ have a join (smallest common 
upper bound) and a meet (largest common lower bound), which will
be denoted by $\pi \vee \rho$ and $\pi \wedge \rho$, respectively.
\end{notation}

\begin{remark}    \label{rem:3.2}
A block $W$ of a partition $\pi \in NC(n)$ is called an 
{\em interval-block} if it is of the form $W = [p,q] \cap \bZ$
for some $1 \leq p \leq q \leq n$. Every non-crossing partition 
has interval-blocks, and it is actually easy to check that the 
following more refined statement holds: let $\pi$ be in 
$NC(n)$, let $V$ be a block of $\pi$, and let $i < j$ be two 
elements of $V$ which are consecutive in $V$ (in the sense that 
$(i,j) \cap V \neq  \emptyset$). If $j \neq i+1$ (hence the 
interval $(i,j)$ contains some integers) then there exists an
interval-block $W$ of $\pi$ such that $W \subseteq (i,j)$.
\end{remark}

\begin{notation}    \label{def:3.3}
The lattice of non-crossing partitions of type B of $2n$ elements
will be denoted by $\ncb (n)$. Following the paper of Reiner 
\cite{R1997} where $\ncb (n)$ was introduced, it is customary to
denote the $2n$ elements that are being partitioned as 
$1, \ldots , n$ and $-1, \ldots , -n$, taken in the order 
$1 < \cdots < n < -1 < \cdots < -n$. If we denote by 
$NC( \pm n )$
the lattice
\footnote{ $NC( \pm n )$ is thus just a copy of $NC(2n)$, where 
one puts different labels on some of the $2n$ points that are 
being partitioned.}
of all non-crossing partitions of the ordered set
$\{ 1, \ldots , n \} \cup \{ -1 , \ldots , -n \}$, then $\ncb (n)$
consists of those partitions $\tau \in NC( \pm n )$ which have the
symmetry property that 
\[
\mbox{ $\bigl( V$ is a block of $\tau \bigr)$ }
\Rightarrow
\mbox{ $\bigl( -V$ is a block of $\tau \bigr)$ }
\]
(with $-V := \{ -v \mid v \in V \} \subseteq 
\{ 1, \ldots , n \} \cup \{ -1 , \ldots , -n \}$). $\ncb (n)$ 
inherits from $NC( \pm n )$ the partial order by reverse refinement,
and is closed under the operations $\vee, \wedge$, hence is a 
sublattice of $NC ( \pm n )$. Note also that $\ncb (n)$ contains 
the minimal and maximal elements of $\ncb (n)$, which will be denoted
as $0_{\pm n}$ and $1_{\pm n}$, respectively.

A block $Z$ of a partition $\tau \in \ncb (n)$ is called a 
{\em zero-block} when it satisfies the condition $Z = -Z$. The set
$\{ \tau \in \ncb (n) \mid \tau$ has zero-blocks$\}$ will be 
denoted by $\nczb (n)$. Due to the non-crossing property, it is 
immediate that every $\tau \in \nczb (n)$ has exactly one 
zero-block (hence it is justified to talk about ``the zero-block''
of $\tau$).
\end{notation}

\begin{remark}    \label{rem:3.4}
{\em (Kreweras complementation.)} An important ingredient in the 
study of the lattice $NC(n)$ is a special anti-automorphism 
$\Kr : NC(n) \to NC(n)$, called the {\em Kreweras complementation 
map} (see pp. 147-148 in \cite{NS2006}). Since 
$NC( \pm n ) \simeq NC(2n)$, one also has such a map $\Kr$ on
$NC( \pm n )$. (All occurrences of Kreweras complementation 
maps in this paper will be denoted in the same way, by ``$\Kr$''.) 
Moreover, the sublattice $\ncb (n) \subseteq NC( \pm n )$ turns 
out to be invariant under the $\Kr$ map of $NC( \pm n)$, hence one 
can talk about the Kreweras complementation map on $\ncb (n)$ as 
well. It is easily checked that $\Kr : \ncb (n) \to \ncb (n)$ maps 
the sets $\nczb (n)$ and $\ncb (n) \setminus \nczb (n)$ 
bijectively onto each other (see Section 1.2 of \cite{BGN2003}).
\end{remark}

\begin{remark}    \label{rem:3.5}
{\em (Absolute value map.)} Let 
$\abs : \{ 1, \ldots , n \} \cup \{ -1 , \ldots , -n \} \to
\{ 1, \ldots , n \}$ denote the absolute value map sending $\pm i$
to $i$ for $1 \leq i \leq n$. In \cite{BGN2003} it was observed 
that it makes sense to extend the concept of ``absolute value'' 
to non-crossing partitions. That is, for $\tau \in \ncb (n)$
it makes sense to define $\abs ( \tau ) \in NC(n)$
to be the partition of $\{ 1, \ldots , n \}$ into blocks of the 
form $\abs (V)$, $V \in \tau$. Moreover, Section 1.4 of 
\cite{BGN2003} puts into evidence the remarkable fact that the
map $\abs : \ncb (n) \to NC(n)$ so defined is an $(n+1)$-to-$1$ 
map, and explains precisely how to find the $n+1$ partitions in  
$\abs^{-1} ( \pi )$, for a given $\pi \in NC(n)$. A part of 
this result which is important for the present paper is that for 
every $\pi \in NC(n)$ and $V \in \pi$ there exists a unique 
$\tau \in \nczb (n)$ such that $\abs ( \tau ) = \pi$ and such 
that the zero-block $Z$ of $\tau$ has $\abs (Z) = V$.
Clearly, this can be
rephrased by saying that we have a bijection
\begin{equation}   \label{eqn:3.51}
\left\{  \begin{array}{rcl}
\nczb (n)  & \longrightarrow  & 
             \{ ( \pi , V ) \mid \pi \in NC(n), \ V \in \pi \}  \\
\tau       & \mapsto          & 
             \bigl( \abs ( \tau ), \abs (Z) \bigr)              \\
           &                  &
         \mbox{ (where $Z :=$ the unique zero-block of $\tau$).}
\end{array}  \right.
\end{equation}
Moreover, for every $\pi \in NC(n)$, the $n+1 - | \pi |$ 
partitions in $\abs^{-1} ( \pi )$ that are not accounted by 
(\ref{eqn:3.51}) are all from $\ncb (n) \setminus \nczb (n)$, and 
are naturally indexed by the blocks of $\Kr ( \pi )$. For the 
explanation of why (and how) this happens, we refer to the 
Remark on p. 2270 of \cite{BGN2003}. 
\end{remark}

\begin{remark}    \label{rem:3.6}
{\em (M\"obius functions.)}
We will use the notation ``$\moebA$'' for the M\"obius functions of 
the lattices $NC(n)$. The value $\moebA ( \pi , \rho )$ for 
$\pi \leq \rho$ in $NC(n)$ can be given explicitly, as a product of
signed Catalan numbers (see p. 163 in \cite{NS2006}). In the 
present paper we will not need the concrete values 
$\moebA ( \pi , \rho )$, but only the M\"obius inversion formula; 
this says that if we have two families of vectors 
$\{ f_{\pi} \mid \pi \in NC(n) \}$ and
$\{ g_{\pi} \mid \pi \in NC(n) \}$ in the same vector space over
$\bC$, then the relations 
\begin{equation}   \label{eqn:3.61}
g_{\rho} = \sum_{ \pi \in NC(n), \ \pi \leq \rho } \ \ 
f_{\pi}, \ \ \ \forall \, \rho \in NC(n)
\end{equation}
are equivalent to
\begin{equation}   \label{eqn:3.62}
f_{\rho} = \sum_{ \pi \in NC(n), \ \pi \leq \rho } \ \ 
\moebA ( \pi, \rho ) \cdot g_{\pi}, 
\ \ \ \forall \, \rho \in NC(n).
\end{equation}

We will use the notation ``$\moebB$'' for the M\"obius functions of 
the lattices $\ncb (n)$. The explicit values 
$\moebB ( \sigma , \tau )$ for $\sigma \leq \tau$ in $\ncb (n)$ can 
be read off from the considerations in Section 3 of \cite{R1997}. 
Here we will only need a simple connection between the types 
A and B, saying that 
\begin{equation}   \label{eqn:3.63}
\Bigl( \sigma \leq \tau \mbox{ in } \nczb (n) \Bigr) \
\Rightarrow \
\moebB ( \sigma , \tau ) =
\moebA \bigl( \abs ( \sigma ) , \abs ( \tau ) \bigr) .
\end{equation}
For the proof of (\ref{eqn:3.63}) one observes that $\abs$ gives a 
poset isomorphism between the intervals
$[ \sigma , \tau ] \subseteq \ncb (n)$ and 
$[ \abs ( \sigma ) , \abs ( \tau ) ] \subseteq NC(n)$, then uses
the fact that the values $\moebB ( \sigma , \tau )$ and 
$\moebA \bigl( \abs ( \sigma ) , \abs ( \tau ) \bigr)$ only depend 
on the isomorphism classes (in the category of posets) of these
intervals.
\end{remark}

$\ $

\begin{center}
{\bf 3.2 Non-crossing cumulants, in the usual $\bC$-valued setting}
\end{center}

The following notation for ``restrictions of $n$-tuples'' will be
used throughout the whole paper.

\begin{notation}   \label{def:3.11}
Let $(a_1, \ldots , a_n)$ be an $n$-tuple of elements in a set
$\cA$, and let $V = \{ v_1, \ldots , v_m \}$ be a non-empty subset
of $\{ 1, \ldots , n \}$, with $v_1 < \cdots < v_m$. Then we 
denote
\begin{equation}    \label{eqn:3.111}
(a_1, \ldots , a_n) \mid V  :=
( a_{v_1}, \ldots , a_{v_m} ) \in \cA^m .
\end{equation}
\end{notation}

\begin{definition}   \label{def:3.12}
Let $( \cA , \varphi )$ be a noncommutative probability space.
The multilinear functionals 
$( \, \kk_n : \cA^n \to \bC \, )_{n \geq 1}$ defined by
\begin{equation}   \label{eqn:3.121}
\begin{array}{lr}
\kappa_n (a_1, \ldots , a_n ) = &
\sum_{ \pi \in NC(n) } 
\ \Bigl( \, \moebA ( \pi , 1_n ) \cdot \prod_{V \in \pi} 
\varphi_{ |V| } ( \, (a_1, \ldots , a_n) \mid V \, ) \, \Bigr),  \\
    & \mbox{for $n \geq 1$ and $a_1, \ldots , a_n \in \cA$}
\end{array}
\end{equation}	
are called the {\em non-crossing cumulant functionals} associated 
to $( \cA , \varphi )$.
\end{definition}

The importance of non-crossing cumulants for free probability 
theory comes from the following theorem, originally found in 
\cite{S1994} (see also the detailed presentation in Lecture 11
of \cite{NS2006}).

\begin{theorem}   \label{thm:3.13}
Let $( \cA , \varphi )$ be a noncommutative probability space and 
let $\cA_1, \ldots , \cA_k$ be unital subalgebras of $\cA$. The 
following statements are equivalent:

\noindent
(1) $\cA_1, \ldots , \cA_k$ are free.

\noindent
(2) For every $n \geq 2$, for every 
$i_1, \ldots , i_n \in \{ 1, \ldots , k \}$ which are not 
all equal to each other, and for every
$a_1 \in \cA_{i_1}, \ldots , a_n \in \cA_{i_n}$, one has that 
$\kappa_n (a_1, \ldots , a_n) = 0$.
\end{theorem}

\begin{remark}      \label{rem:3.14}
Let $( \cA , \varphi )$ be a noncommutative probability space. For
every $n \geq 1$ let us consider the multiplication map
\begin{equation}   \label{eqn:3.141}
\mult_n : \cA^n \to \cA, \ \ \
\mult_n (a_1, \ldots , a_n) = a_1 \cdots a_n,
\end{equation}
and let us denote $\varphi_n := \varphi \circ \mult_n$. 
The multilinear functionals $( \varphi_n )_{n \geq 1}$ are called
the {\em moment functionals} of $( \cA , \varphi )$. For every
$n \geq 1$ and $\pi \in NC(n)$ let us next define a multilinear
functional $\phiA_{\pi} : \cA^n \to \bC$ by
\footnote{ The superscript ``$(A)$'' is used in anticipation of 
the fact that some multilinear functionals $\phiB_{\tau}$ with
$\tau \in \ncb (n)$ will appear in Section 6 of the paper. }
\begin{equation}   \label{eqn:3.142}
\phiA_{\pi} ( a_1, \ldots , a_n ) :=
\prod_{V \in \pi} \varphi_{ |V| } 
 \bigl( \, (a_1, \ldots , a_n) \mid V \, \bigr),
\ \ \ a_1, \ldots , a_n \in \cA .
\end{equation}
Then Definition \ref{def:3.12} can be rephrased as saying that
\begin{equation}   \label{eqn:3.143}
\kappa_n = \sum_{\pi \in NC(n)} \ \moebA ( \pi , 1_n ) \cdot
\phiA_{\pi} \in \fM_n,
\end{equation}
where $\fM_n$ denotes the vector space of multilinear functionals
from $\cA^n$ to $\bC$. Moreover, if for every $\pi \in NC(n)$ we 
introduce (by analogy with (\ref{eqn:3.142})) a functional 
$\kappaA_{\pi} \in \fM_n$ defined by
\begin{equation}   \label{eqn:3.144}
\kappaA_{\pi} ( a_1, \ldots , a_n ) :=
\prod_{V \in \pi} \kappa_{ |V| } 
 \bigl( \, (a_1, \ldots , a_n) \mid V \, \bigr),
\ \ \ a_1, \ldots , a_n \in \cA ,
\end{equation}
then it is not hard to see that the formula (\ref{eqn:3.143}) for 
$\kappa_n$ extends to
\begin{equation}   \label{eqn:3.145}
\kappa_{\rho} = \sum_{\pi \in NC(n), \ \pi \leq \rho} \ \ 
\moebA ( \pi , \rho ) \cdot \phiA_{\pi}, \ \ \ 
\forall \, \rho \in NC(n).
\end{equation}
Thus for a given $n \geq 1$, the families of functionals 
$\{ \kappaA_{\pi} \mid \pi \in NC(n) \}$ and 
$\{ \phiA_{\pi} \mid \pi \in NC(n) \}$ are exactly as in the above
Remark \ref{rem:3.6}. Equation (\ref{eqn:3.145}) and its 
equivalent counterpart which express $\phiA_{\rho}$ as the sum 
of the functionals $\{ \kappaA_{\pi} \mid \pi \leq \rho \}$ go 
under the name of {\em non-crossing moment-cumulant formulas}
for $( \cA , \varphi )$.
\end{remark}

$\ $

\begin{center}
{\bf 3.3 Non-crossing cumulants in the $\bG$-valued setting}
\end{center}

\begin{remark}     \label{rem:3.21}
We will work with the Grassman algebra $\bG$ from subsection 1.1, 
and with the maps $\bo , \so : \bG \to \bC$ defined in subsection 
1.2. It is immediate that the multiplication of $\bG$ is 
commutative, and that the ``body'' map $\bo : \bG \to \bC$ is a 
homomorphism of unital algebras. Concerning how the ``soul'' map 
$\so$ behaves with respect to multiplication, we record the 
immediate formula 
\begin{equation}     \label{eqn:3.211}
\so ( \gamma_1 \cdots \gamma_n ) =
\sum_{i=1}^n \Bigl( \, \so ( \gamma_i ) \cdot 
\prod_{ \begin{array}{c}
{\scriptstyle 1 \leq j \leq n, } \\
{\scriptstyle j \neq i}
\end{array} } \, \bo ( \gamma_j ) \, \Bigr) , \ \ \forall 
\, n \geq 1, \ \forall \, \gamma_1, \ldots ,\gamma_n \in \bG .
\end{equation}
\end{remark}

\begin{notation}   \label{def:3.22}
{\em For the rest of this subsection we fix a pair
$( \cA , \phitild )$} where $\cA$ is a unital algebra over $\bC$
and $\phitild : \cA \to \bG$ is $\bC$-linear with
$\phitild ( 1_{\cA} ) = 1$. 
In connection to this $\phitild$ we will repeat all the 
constructions of functionals described in Remark \ref{rem:3.14}, 
with the only difference that the range space of these 
functionals is now $\bG$. So for every $n \geq 1$ we put 
$\phitild_n = \phitild \circ \mult_n : \cA^n \to \bG$,  
where $\mult_n : \cA^n \to \cA$ is the same as in 
Equation (\ref{eqn:3.141}). Then for every $\pi \in NC(n)$ 
we define $\phitild_{\pi} : \cA^n \to \bG$ by
\begin{equation}   \label{eqn:3.221}
\phitild_{\pi} ( a_1, \ldots , a_n ) :=
\prod_{V \in \pi} \phitild_{ |V| } 
 \bigl( \, (a_1, \ldots , a_n) \mid V \, \bigr),
\ \ \ a_1, \ldots , a_n \in \cA .
\end{equation}
This is followed by defining a family of cumulant functionals 
$( \, \ktild_n : \cA^n \to \bG \, )_{n \geq 1}$, where
\begin{equation}   \label{eqn:3.222}
\ktild_n = \sum_{ \pi \in NC(n) } 
\, \moebA ( \pi , 1_n ) \cdot \phitild_{\pi} , \ \ n \geq 1.
\end{equation}	
Finally, for every $\pi \in NC(n)$ we define 
$\ktild_{\pi} : \cA^n \to \bG$ by
\begin{equation}   \label{eqn:3.223}
\ktild_{\pi} ( a_1, \ldots , a_n ) :=
\prod_{V \in \pi} \ktild_{ |V| } 
 \bigl( \, (a_1, \ldots , a_n) \mid V \, \bigr),
\ \ \ a_1, \ldots , a_n \in \cA .
\end{equation}
It is easily seen that, exactly as in the $\bC$-valued case 
from Remark \ref{rem:3.14}, the families of functionals 
$\{ \ktild_{\pi} \mid \pi \in NC(n) \}$ and 
$\{ \phitild_{\pi} \mid \pi \in NC(n) \}$ are related by 
moment-cumulant formulas (i.e. by summation formulas as shown 
in Equations (\ref{eqn:3.61}), (\ref{eqn:3.62}) of
Remark \ref{rem:3.6}). We only record here the special case of
moment-cumulant formula which expresses $\phitild_{1_n}$ as 
a sum of cumulant functionals, and thus says that
\begin{equation}    \label{eqn:3.224}
\phitild (a_1 \cdots a_n ) =
\sum_{ \pi \in NC(n) } 
\ \ktild_{\pi} (a_1, \ldots , a_n)  \in \bG, \ \ \
\forall \, a_1, \ldots , a_n \in \cA .
\end{equation}	
\end{notation}

\begin{remark}    \label{rem:3.23}
A natural question concerning $( \cA , \phitild )$ is whether
the analogue of Theorem \ref{thm:3.13} is holding in this 
framework. As will be explained in detail in Remark \ref{rem:4.8} 
below, both conditions $(1)$ and $(2)$ from the statement of 
Theorem \ref{thm:3.13} can be faithfully transcribed in the 
context of $( \cA , \phitild )$, but then they are no longer 
equivalent to each other -- the implication 
$(2) \Rightarrow (1)$ still holds, but its converse does not. 

In the remaining part of this subsection we will point out 
two other facts from the theory of usual non-crossing cumulants 
where (unlike for Theorem \ref{thm:3.13}) both the statement 
and the proof can be transcribed without any problems from 
usual $\bC$-valued framework to the $\bG$-valued framework
of $( \cA , \phitild )$.
\end{remark}

\begin{proposition}    \label{prop:3.24}
One has that $\ktild_n (a_1, \ldots , a_n) = 0$
whenever $n \geq 2$, $a_1, \ldots , a_n \in \cA$,
and there exists $1 \leq m \leq n$ such that $a_m \in \bC 1_{\cA}$. 
\end{proposition}

\begin{proof}
This is the analogue of Proposition 11.15 in \cite{NS2006}. It is 
straightforward (left to the reader) to see that the proof shown 
on p. 182 of \cite{NS2006} goes without any changes to the 
$\bG$-valued framework.
\end{proof}

\begin{proposition}    \label{prop:3.25}
Let $x_1, \ldots , x_s$ be in $\cA$ and consider some products 
of the form
\[
a_1 = x_1 \cdots x_{s_1}, \ a_2 = x_{s_1 +1} \cdots x_{s_2}, 
\ \ldots , \ a_n = x_{s_{n-1} +1} \cdots x_{s_n},
\]
where $1 \leq s_1 < s_2 < \cdots < s_n = s$. Then 
\begin{equation}    \label{eqn:3.251}
\ktild_n (a_1, \ldots , a_n) = 
\sum_{ \begin{array}{c}
{\scriptstyle \pi \in NC(s) \ such} \\
{\scriptstyle that \ \pi \vee \theta = 1_s}
\end{array}  } \ \ \ktild_{\pi} (x_1, \ldots , x_s),
\end{equation}
where $\theta \in NC(s)$ is the partition with interval blocks
$\{ 1, \ldots , s_1 \}$,
$\{ s_1 + 1, \ldots , s_2 \} , \ldots ,
\{ s_{n-1} + 1, \ldots , s_n \}$.
\end{proposition}

\begin{proof} This is the analogue of Theorem 11.20 in 
\cite{NS2006}, and the proof of this theorem (as shown on pp.
178-180 of \cite{NS2006}) goes without any changes to the 
$\bG$-valued framework.
\end{proof}

$\ $

$\ $

\begin{center}
{\bf\large 4. Infinitesimal cumulants and the proof of 
Theorem \ref{thm:1.2} }
\end{center}
\setcounter{section}{4}
\setcounter{equation}{0}
\setcounter{theorem}{0}

\begin{notation}   \label{def:4.1}
{\em Throughout this whole section we fix an incps 
$( \cA , \varphi , \phiprim )$}. We will use the notation 
``$\kappa_n$'' for the non-crossing cumulant functionals
associated to $\varphi$, as described in Section 3.2. 
Moreover, we will denote, same as in the introduction:
\[
\phitild = \varphi + \ee \phiprim : \cA \to \bG
\]
and we will consider the family of non-crossing cumulant 
functionals $( \ktild_n : \cA^n \to \bG )_{n \geq 1}$ which 
are associated to $\phitild$ as in Section 3.3.
\end{notation}

\begin{definition}   \label{def:4.2}
For every $n \geq 1$, consider the multilinear functional
$\kappa_n ' : \cA^n \to \bC$ defined by the formula
\begin{equation}    \label{eqn:4.21}
\kappa_n ' (a_1, \ldots , a_n) = 
\end{equation}
\[
\sum_{\pi \in NC(n)} 
\sum_{V \in \pi} \Bigl[ \ \moeb ( \pi , 1_n ) \
\phiprim_{|V|} ( \, (a_1, \ldots , a_n) \mid V \, ) \cdot
\prod_{ \begin{array}{c}
{\scriptstyle W \in \pi}  \\
{\scriptstyle W \neq V}
\end{array} } \ 
\varphi_{|W|} ( \, (a_1, \ldots , a_n) \mid W \, ) \ \Bigr] ,
\]
for $a_1, \ldots , a_n \in \cA$. The functionals $\kappa_n '$ 
will be called {\em infinitesimal non-crossing cumulant functionals} 
associated to $( \cA , \varphi , \phiprim )$.
\end{definition}

A moment's thought shows that Equation (\ref{eqn:4.21}) is indeed
obtained from the fomula (\ref{eqn:3.121}) defining $\kappa_n$, 
where one uses the formal derivation procedure announced in 
subsection 1.2 of the introduction. 

We next make precise (in Propositions \ref{prop:4.20},
\ref{prop:4.4} and Remark \ref{rem:4.3}) the equivalence between
Definition \ref{def:4.2} and the other facets of $\kappa_n '$ that 
were mentioned in subsection 1.2.

\begin{proposition}         \label{prop:4.20}
Suppose that $\varphi , \phiprim$ are the infinitesimal limit 
of a family $\{ \varphi_t \mid t \in T \}$, in the sense described 
in Equation (\ref{eqn:1.23}). Let us use the notation $\kk^{(t)}_n$ 
for the non-crossing cumulant functional of $\varphi_{t}$, for 
$t \in T$ and $n \geq 1$. Then for every 
$n \geq 1$ and every $a_1, \ldots , a_n \in \cA$ one has that
\[
\kk_n (a_1, \ldots , a_n) =
\lim_{t \to 0} \kk^{(t)}_n (a_1, \ldots , a_n),
\]
and
\[
\kappa_n ' (a_1, \ldots , a_n) =
\Bigl[ \ \frac{d}{dt} \kk^{(t)}_n (a_1, \ldots , a_n) \, \Bigr]
\ \vline \ { }_{ { }_{t=0} }.
\]
\end{proposition}

\begin{proof} Fix $n \geq 1$ and $a_1, \ldots , a_n \in \cA$. 
For every $t \in T$ we have that 
\begin{equation}        \label{eqn:4.201}
\kk^{(t)}_n (a_1, \ldots , a_n) = \sum_{\pi \in NC(n)} \,
\moebA ( \pi , 1_n ) \cdot \prod_{V \in \pi}
\varphi_t \bigl( \, (a_1, \ldots , a_n) \mid V \, \bigr) .
\end{equation}
From (\ref{eqn:4.201}) it is clear that 
$\lim_{t \to 0} \kk^{(t)}_n (a_1, \ldots , a_n) = 
\kk_n (a_1, \ldots , a_n)$. Moreover, it is immediate that the 
function of $t$ appearing on the right-hand side of 
(\ref{eqn:4.201}) has a derivative at $0$; and upon using linearity 
and the Leibnitz formula to compute this derivative, one 
obtains precisely the formula (\ref{eqn:4.21}) that 
defined $\kk_n ' (a_1, \ldots , a_n)$.
\end{proof}

\begin{remark}    \label{rem:4.3}
As observed in Remark \ref{rem:3.5}, the set 
$\{ ( \pi , V ) \mid \pi \in NC(n), V \in \pi \}$
which indexes the sum on the right-hand side of Equation 
(\ref{eqn:4.21}) is the image of $\nczb (n)$ via the bijection
$\bigl( \tau \in \nczb (n) \mbox{ with zero-block } Z \bigr) 
\mapsto ( \, \abs ( \tau ), \abs (Z) \, )$. 
When $\tau$ and $( \pi , V )$ correspond to each other via 
this bijection, we have that $\moebB ( \tau , 1_{\pm n} )$ = 
$\moebA ( \pi , 1_n )$ (cf. implication (\ref{eqn:3.63}) in 
Remark \ref{rem:3.6}); moreover, the rest of 
the product indexed by $( \pi , V )$ on the right-hand side of 
Equation (\ref{eqn:4.21}) is precisely equal to 
$\phiB_{\tau} ( a_1, \ldots , a_n)$, where we anticipate here 
the notation $\phiB_{\tau}$ from Equation (\ref{eqn:6.22}).
In conclusion, the change of variable from $(V , \pi )$ to 
$\tau$ converts (\ref{eqn:4.21}) into a summation formula 
``of type B'',
\begin{equation}    \label{eqn:4.33}
\kappa_n ' = \sum_{\tau \in \nczb (n)} \
\moebB ( \tau , 1_{\pm n} ) \cdot \phiB_{\tau} .
\end{equation}
It is easy to see that (\ref{eqn:4.33}) is equivalent to a 
plain summation formula which writes $\phiprim (a_1 \cdots a_n)$ 
in terms of cumulants (cf. Remark \ref{rem:6.5} below, where 
one also sees that the absence of terms indexed by partitions 
from $\ncb (n) \setminus \nczb (n)$ is caused by the fact that 
$\phiprim ( 1_{\cA} ) = 0$).
\end{remark}

\begin{proposition}   \label{prop:4.4}
For every $n \geq 1$ one has that 
$\bo \ \ktild_n  = \kappa_n$ and
$\so \ \ktild_n  = \kappa_n '$.
\end{proposition}

\begin{proof}
For the first statement we only have to take the body part on 
both sides of Equation (\ref{eqn:3.222}) and use the fact that 
$\bo : \bG \to \bC$ is a homomorphism of unital algebras. 
For the second statement we take soul parts in (\ref{eqn:3.222}) 
and then use the multiplication formula (\ref{eqn:3.211}).
\end{proof}

We now go to Theorem \ref{thm:1.2}. Note that, in view of 
Proposition \ref{prop:4.4}, the equalities 
``$\kappa_n (a_1, \ldots , a_n) = \kappa_n ' (a_1, \ldots , a_n) 
=0$'' from condition $(2)$ of Theorem \ref{thm:1.2} may be replaced
with ``$\ktild_n (a_1, \ldots , a_n) = 0$''. We will prove Theorem
\ref{thm:1.2} in this alternative form, which is stated below as
Proposition \ref{prop:4.6}.

\begin{lemma}   \label{lemma:4.5}
Suppose that $n$ is a positive integer and $\pi$ is a partition 
in $NC(n)$, such that the following two properties hold:

\noindent
(i) For every $1 \leq i \leq n-1$, the numbers $i$ and $i+1$ do 
not belong to the same block of $\pi$.

\noindent
(ii) $\pi$ has at most one block of cardinality 1.

\noindent
Then $n$ is odd, and $\pi$ is the partition
\[
\Bigl\{ \ \{ 1,n \} , \, \{ 2,n-1 \} , \ldots , 
\{ (n-1)/2, (n+3)/2 \} , \, \{ (n+1)/2 \} \ \Bigr\} .
\]
\end{lemma}

\begin{proof} We will use the observation about interval-blocks
of non-crossing partitions that was recorded in Remark
\ref{rem:3.2}. Clearly, condition (i) implies that $\pi$ 
cannot have interval-blocks $V$ with $|V| \geq 2$; by also 
taking (ii) into account we thus see that $\pi$ has a unique 
interval-block $V_o$, of the form $V_o = \{ p \}$ for 
some $1 \leq p \leq n$.

Let $V$ be a block of $\pi$, distinct from $V_o$. We claim that
\begin{equation}     \label{eqn:4.51}
\mid V \cap [ 1, p ) \mid \ \leq 1 , \ \
\mid V \cap ( p, n ] \mid \ \leq 1.
\end{equation} 
Indeed, assume for instance that we had 
$\mid V \cap [ 1, p ) \mid \ \geq 2$. Then we could find $i,j \in V$
such that $i < j < p$ and $(i,j) \cap V = \emptyset$. Note that 
$j \neq i+1$, due to condition (i); but then, as observed in Remark 
\ref{rem:3.2}, the partition $\pi$ must have an 
interval-block $W \cap (i,j)$, in contradiction to the fact that 
the unique interval-block of $\pi$ is $V_o$.

For every block $V \neq V_o$ of $\pi$ it then follows that 
$\mid V \cap [ 1, p ) \mid \ = \
\mid V \cap ( p, n ] \mid  \ = 1$.
Indeed, if in (\ref{eqn:4.51}) one of the sets 
$V \cap [ 1, p )$, $V \cap ( p, n ]$
would be empty, then it would follow that $|V| =1$ and hypothesis 
(ii) would be contradicted.

The list of blocks of $\pi$ which are distinct from $V_o$ can thus 
be written in the form
\begin{equation}    \label{eqn:4.53}
\left\{  \begin{array}{l}
V_1 = \{ i_1 , j_1 \}, \ldots , V_m = \{ i_m , j_m \} ,  
\ \ \mbox{ where}                                          \\ 
i_1 < p < j_1, \ldots , i_m < p < j_m, \ \mbox{ and } \
i_1 < i_2 < \cdots < i_m.
\end{array}  \right.
\end{equation}
Observe that in (\ref{eqn:4.53}) we must have 
$j_1 > j_2 > \cdots > j_m$. Indeed, if it was true that 
$j_s < j_t$ for some $1 \leq s < t \leq m$, then it would follow
that $i_s < i_t < p < j_s < j_t$, and the blocks $V_s, V_t$
would cross. Hence we have obtained
$i_1 < \cdots < i_m  < p < j_m < \cdots < j_1$; 
together with (\ref{eqn:4.53}), this implies that $n = 2m+1$ 
and that $\pi$ is precisely the partition indicated in the 
lemma. 
\end{proof}

\begin{proposition}    \label{prop:4.6}
Let $\cA_1, \ldots , \cA_k$ be unital subalgebras of $\cA$. The 
following statements are equivalent:

\noindent
(1) $\cA_1, \ldots , \cA_k$ are infinitesimally free in 
$( \cA , \varphi , \phiprim )$.

\noindent
(2) For every $n \geq 2$, for every 
$i_1, \ldots , i_n \in \{ 1, \ldots , k \}$ which are not all 
equal to each other, and for every 
$a_1 \in \cA_{i_1} , \ldots , a_n \in \cA_{i_n}$, one has that 
$\ktild_n (a_1, \ldots , a_n) = 0$.
\end{proposition}

\begin{proof} 
{\em ``$(1) \Rightarrow (2)$''.} We prove the required statement
about cumulants by induction on $n$.
For the base case $n=2$, consider elements $a_1 \in \cA_{i_1}$ and 
$a_2 \in \cA_{i_2}$, where $i_1 \neq i_2$. By using the 
formulas which define $\kappa_2$ and $\kappa_2 '$ and by 
invoking Equations (\ref{eqn:2.21}) and (\ref{eqn:2.22}) 
from Remark \ref{rem:2.2} we find that 
\[
\left\{   \begin{array}{l}
\kappa_2 (a_1, a_2) = \varphi (a_1 a_2) - 
\varphi (a_1)\varphi (a_2) = 0 \ \ \mbox{ and }           \\ 
\kappa_2 ' (a_1, a_2) = \phiprim (a_1 a_2) - 
\phiprim (a_1) \varphi (a_2) - 
\varphi (a_1) \phiprim (a_2) = 0,
\end{array}   \right.
\]
hence $\ktild_2 (a_1, a_2 ) =  \kappa_2 (a_1, a_2) + \ee
\kappa_2 ' (a_1, a_2) = 0$.

We now prove the induction step: assume that the vanishing 
of mixed cumulants is already proved for $1,2,...,n-1$, where 
$n\geq 3$. We consider elements 
$a_1 \in \cA_{i_1} , \ldots , a_n \in \cA_{i_n}$ where 
not all indices $i_1, \ldots , i_n$ are equal to each other,
and we want to prove that $\ktild_n (a_1, \ldots , a_n) = 0$.
By invoking Proposition \ref{prop:3.24} we may replace every 
$a_m$ with $a_m - \varphi (a_m) 1_{\cA}$, $1 \leq m \leq n$,
and therefore assume without loss of generality that 
$\varphi (a_1) = \cdots = \varphi (a_n) = 0$. Observe that 
this implies 
$\phitild (a_p) \phitild (a_q)
= ( \ee \phiprim (a_p) ) \cdot ( \ee \phiprim (a_q) ) = 0$, 
hence that 
\begin{equation}   \label{eqn:4.61} 
\ktild_2 ( a_p, a_q) = \phitild (a_p a_q)
- \phitild (a_p) \phitild (a_q) = \phitild (a_p a_q),
\ \ \forall \, 1 \leq p<q \leq n.
\end{equation}
Another assumption that can be made without loss of generality
is that $i_m \neq i_{m+1}, \forall \, 1 \leq m < n$. Indeed, if 
there exists $1 \leq m < n$ such that $i_m = i_{m+1}$, then we 
invoke the special case of Proposition \ref{prop:3.25} which 
states that 
\begin{equation}   \label{eqn:4.615}
\ktild_{n-1} ( a_1, \ldots , a_m a_{m+1} , \ldots , a_n )
= \ktild_n ( a_1, \ldots , a_n ) + 
\sum_{ \begin {array}{c}
{\scriptstyle \pi \in NC(n) \,with \, | \pi | =2} \\
{\scriptstyle \pi \, separates \, m \, and \,m+1} 
\end{array} } \
\ktild_{\pi} ( a_1, \ldots , a_n ).
\end{equation}
The induction hypothesis gives us that the left-hand side and
every term in the sum on the right-hand side of Equation 
(\ref{eqn:4.615}) are equal to 0, and it follows that 
$\ktild_n ( a_1, \ldots , a_n )$ must be $0$ as well. 

Hence for the rest of the proof of this induction step we will 
assume that $\varphi (a_1) = \cdots = \varphi (a_n) = 0$ and
that $i_1 \neq i_2, \ldots , i_{n-1} \neq i_n$. This makes 
$a_1 , \ldots , a_n$ be exactly as in Definition \ref{def:1.1}, 
so we get that $\varphi (a_1 \cdots a_n) = 0$ and that 
$\phiprim (a_1 \cdots a_n)$ is as described in Equation 
(\ref{eqn:1.15}). In terms of the functional $\phitild$, we have
\begin{equation}     \label{eqn:4.62}
\phitild (a_1 \cdots a_n) = \ee \phiprim (a_1 \cdots a_n) =
\end{equation}
\[
= \left\{  \begin{array}{l}
\ee \varphi (a_1 \, a_n) \varphi (a_2 \, a_{n-1}) \cdots 
       \varphi ( a_{(n-1)/2} \, a_{(n+3)/2} ) \cdot
       \phiprim ( a_{(n+1)/2} ),                                \\
\mbox{$\ \ $} \ \ \mbox{ if $n$ is odd and $i_1 = i_n , 
      i_2 = i_{n-1}, \ldots , i_{(n-1)/2} = i_{(n+3)/2}$, }     \\
0, \ \mbox{ otherwise.}
\end{array}  \right.
\]

Now let us consider the relation (\ref{eqn:3.224}), written in 
the equivalent form
\begin{equation}   \label{eqn:4.63}
\ktild_n ( a_1, \ldots , a_n ) = \phitild (a_1 \cdots a_n) -
\sum_{ \begin{array}{c}
{\scriptstyle \pi \in NC(n),}  \\
{\scriptstyle \pi \neq 1_n} 
\end{array}  } \ \ktild_{\pi} ( a_1, \ldots , a_n ).
\end{equation}
Observe that if a partition $\pi \in NC(n)$ has two distinct 
blocks $\{ p \}, \{ q \}$ of cardinality one, then the term 
indexed by $\pi$ on the right-hand side of (\ref{eqn:4.63}) 
vanishes, because it contains the subproduct 
$\ktild_1 (a_p) \ktild_1 (a_q) = \phitild (a_p) \phitild (a_q) 
= 0$. On the other hand if $\pi \in NC(n)$ has a block $V$ which 
contains two consecutive numbers $i$ and $i+1$, then the term 
indexed by $\pi$ on the right-hand side of (\ref{eqn:4.63})
vanishes as well, due to the induction hypothesis. Hence the 
sum subtracted on the right-hand side of (\ref{eqn:4.63}) can 
only get non-zero contributions from partitions $\pi \in NC(n)$ 
which satisfy the hypotheses of Lemma \ref{lemma:4.5}; from the 
lemma it then follows that the sum in question is $0$ for $n$ 
even, and is equal to 
\begin{equation}   \label{eqn:4.64}
\ktild_2 (a_1 , a_n) \ktild_2 (a_2 , a_{n-1}) \cdots
\ktild_2 ( a_{(n-1)/2} , a_{(n+3)/2} ) \cdot
\ktild_1 ( a_{(n+1)/2} )
\end{equation}   
for $n$ odd. 

Let us focus for a moment on the quantity that appeared in 
(\ref{eqn:4.64}). The vanishing of mixed cumulants of order 2 
(which is part of our induction hypothesis) implies that this 
quantity vanishes unless $i_1 = i_n $, 
$i_2 = i_{n-1}, \ldots , i_{(n-1)/2} = i_{(n+3)/2}$. In the case 
that the latter equalities of indices hold, we can continue 
(\ref{eqn:4.64}) with 
\[
= \phitild (a_1 a_n) \phitild (a_2 a_{n-1})
\cdots \phitild ( a_{(n-1)/2} a_{(n+3)/2} )
\cdot \phitild ( a_{(n+1)/2} ) \ \ 
\mbox{ (due to (\ref{eqn:4.61})) }
\]
\begin{equation}    \label{eqn:4.65}
= \ee \varphi (a_1 a_n) \varphi (a_2 a_{n-1})
\cdots \varphi ( a_{(n-1)/2} a_{(n+3)/2} )
\cdot \phiprim ( a_{(n+1)/2} ).
\end{equation} 
(The equality (\ref{eqn:4.65}) holds because 
$\phitild ( a_{(n+1)/2} ) = \ee \phiprim ( a_{(n+1)/2} )$, and 
due to how the multiplication on $\bG$ works.)

So all in all, what we have obtained is that 
\begin{equation}   \label{eqn:4.66}
\ktild_n ( a_1, \ldots , a_n ) = 
\end{equation}
\[
= \left\{   \begin{array}{l}
\phitild (a_1 \cdots a_n) - \ee \varphi (a_1 a_n) 
   \varphi (a_2 a_{n-1}) \cdots \varphi ( a_{(n-1)/2} a_{(n+3)/2} )
   \cdot \phiprim ( a_{(n+1)/2} ),                           \\
\mbox{$\ \ $ if $n$ is odd and $i_1 = i_n , 
        i_2 = i_{n-1}, \ldots , i_{(n-1)/2} = i_{(n+3)/2}$,}  \\
\phitild (a_1 \cdots a_n), \ \ \mbox{  otherwise.}
\end{array}  \right.
\]
By comparing Equations (\ref{eqn:4.66}) and (\ref{eqn:4.62})
we see that, in all cases, we have 
$\ktild_n ( a_1, \ldots , a_n ) = 0$.
This concludes the induction argument, and the proof of the 
implication $(1) \Rightarrow (2)$ of the proposition.

\vspace{6pt}

{\em ``$(2) \Rightarrow (1)$''.} Consider indices 
$i_1, \ldots , i_n \in \{ 1, \ldots , k \}$ and elements 
$a_1 \in \cA_{i_1} , \ldots , a_n \in \cA_{i_n}$ such that 
$i_1 \neq i_2, \ldots , i_{n-1} \neq i_n$ and such that 
$\varphi (a_1) = \cdots = \varphi (a_n) = 0$. We have to prove 
that $\varphi (a_1 \cdots a_n) = 0$ and that 
$\phiprim (a_1 \cdots a_n)$ is as described in formula 
(\ref{eqn:1.15}) from Definition \ref{def:1.1}. To this end 
we consider the $\bG$-valued moment $\phitild (a_1 \cdots a_n)$ =
$\varphi (a_1 \cdots a_n) + \ee \phiprim (a_1 \cdots a_n)$,
and write it in terms of $\bG$-valued cumulants as in 
subsection 3.3:
\begin{equation}    \label{eqn:4.67}
\phitild (a_1 \cdots a_n) = \sum_{\pi \in NC(n)} \
\prod_{V \in \pi} \
\ktild_{|V|} ( \, (a_1, \ldots , a_n) \mid V ).
\end{equation}
An argument very similar to the one used in the proof of the 
implication $(1) \Rightarrow (2)$ above shows that the sum on 
the right-hand side of (\ref{eqn:4.67}) can only get non-zero
contributions from partitions $\pi \in NC(n)$ which satisfy 
the hypotheses of Lemma \ref{lemma:4.5}. If $n$ is even then 
there is no such partition, and we obtain 
$\phitild (a_1 \cdots a_n) =0$. If $n$ is odd, then the sum in 
(\ref{eqn:4.67}) reduces to only one term and we obtain that 
\begin{equation}    \label{eqn:4.68}
\phitild (a_1 \cdots a_n) =
\ktild_2 (a_1 , a_n) \ktild_2 (a_2 , a_{n-1}) \cdots
\ktild_2 ( a_{(n-1)/2} , a_{(n+3)/2} ) \cdot
\ktild_1 ( a_{(n+1)/2} ).
\end{equation}
Moreover, in the case when $n$ is odd, the hypothesis that mixed 
cumulants vanish gives us that the right-hand side of 
(\ref{eqn:4.68}) is equal to $0$ unless we have 
$i_1 = i_n , \ldots , i_{(n-1)/2} = i_{(n+3)/2}$. And finally, 
if the latter equalities of indices hold, then the
right-hand side of (\ref{eqn:4.68}) gets converted into 
$\ee \varphi (a_1 a_n) \varphi (a_2 a_{n-1})
\cdots \varphi ( a_{(n-1)/2} a_{(n+3)/2} )
\cdot \phiprim ( a_{(n+1)/2} )$,
by the same argument that led to (\ref{eqn:4.65}) in the proof of 
the implication $(1) \Rightarrow (2)$. The conclusion is that
$\varphi (a_1 \cdots a_n) = 0$ (in all cases), and that 
$\phiprim (a_1 \cdots a_n)$ is as in Equation (\ref{eqn:1.15}), as 
required.
\end{proof}

\begin{corollary}    \label{cor:4.7}
Let $\cX_1, \ldots , \cX_k$ be subsets of $\cA$. The 
following statements are equivalent:

\noindent
(1) $\cX_1, \ldots , \cX_k$ are infinitesimally free in 
$( \cA , \varphi , \phiprim )$.

\noindent
(2) For every $n \geq 2$, for every 
$i_1, \ldots , i_n \in \{ 1, \ldots , k \}$ which are not all 
equal to each other, and for every 
$x_1 \in \cX_{i_1} , \ldots , x_n \in \cX_{i_n}$, one has that 
$\ktild_n (x_1, \ldots , x_n) = 0$.
\end{corollary}

\begin{proof} This is a faithful copy of the proof giving the 
analogous result over $\bC$ (cf. Theorem 11.20 in \cite{NS2006}). 
For the reader's convenience, we repeat here the highlights of 
the argument. Let $\cA_i$ denote the unital subalgebra of $\cA$ 
generated by $\cX_i$, $1 \leq i \leq k$. The infinitesimal freeness 
of $\cX_1 , \ldots , \cX_k$ is by definition equivalent to the 
one of $\cA_1, \ldots , \cA_k$, hence to the fact that condition 
(2) from Proposition \ref{prop:4.6} holds. We must thus prove that 
``(2) in Proposition \ref{prop:4.6}'' is equivalent to
``(2) in Corollary \ref{cor:4.7}''. The implication ``$\Rightarrow$'' 
is trivial. For ``$\Leftarrow$'' it suffices (by multilinearity of 
$\ktild_n$ and Proposition \ref{prop:3.24}) to prove that 
$\ktild_n (a_1, \ldots , a_n) = 0$ when 
\begin{equation}    \label{eqn:4.71}
a_1 = x_1 \cdots x_{s_1}, \
a_2 = x_{s_1 +1} \cdots x_{s_2}, \  \ldots , \ a_n = 
x_{s_{n-1} +1} \cdots x_{s_n}
\end{equation}
for $n \geq 2$ and $1 \leq s_1 < s_2 < \cdots < s_n$, 
where $x_1, \ldots , x_{s_1} \in \cX_{i_1}, \
x_{s_1 +1}, \ldots , x_{s_2} \in \cX_{i_2}, \ldots ,
x_{s_{n-1} +1}, \ldots, x_{s_n} \in \cX_{i_n}$,
and where the indices $i_1, \ldots , i_n$ are not all equal to 
each other. But for $a_1, \ldots , a_n$ as in (\ref{eqn:4.71}), 
Proposition \ref{prop:3.25} gives us the cumulant
$\ktild_n (a_1, \ldots , a_n)$ as a sum of cumulants
$\ktild_{\pi} (x_1, \ldots , x_{s_n})$; and a direct
combinatorial analysis (exactly as on p. 186 of \cite{NS2006})
shows that all the latter cumulants vanish because of 
condition (2) form Corollary \ref{cor:4.7}.
\end{proof}

\begin{remark}    \label{rem:4.8} 
Since the functional $\phitild : \cA \to \bG$ and its associated 
cumulants $\ktild_n$ play such a central role in the proof of 
Theorem \ref{thm:1.2}, it is natural to ask: can't one actually 
characterize infinitesimal 
freeness by the same kind of moment condition as in the definition 
of usual freeness, with the only modification that one now uses
$\phitild$ instead of $\varphi$? To be precise, consider the 
following condition which a family of unital subalgebras 
$\cA_1, \ldots , \cA_k \subseteq \cA$ may or may not satisfy:
\begin{equation}   \label{eqn:4.81}
\left\{  \begin{array}{l}
\mbox{For every $n \geq 1$ and $1 \leq i_1, \ldots , i_n \leq k$
      such that $i_1 \neq i_2, \ldots , i_{n-1} \neq i_n$,}      \\
\mbox{and every $a_1 \in \cA_{i_1}, \ldots , a_n \in \cA_{i_n}$ 
    such that $\phitild (a_1) = \cdots = \phitild (a_n) = 0$,}   \\
\mbox{one has that $\phitild (a_1 \cdots a_n)= 0$.}
\end{array}  \right.
\end{equation}
Isn't then condition (\ref{eqn:4.81}) equivalent to infinitesimal 
freeness?  

On the positive side it is immediate, directly from Definition 
\ref{def:1.1}, that (\ref{eqn:4.81}) is indeed implied by 
infinitesimal freeness. However, the converse statement is not 
true: it may happen that (\ref{eqn:4.81}) is satisfied and yet 
$\cA_1, \ldots , \cA_k$ are not infinitesimally free. What causes 
this to happen is that one cannot generally ``center'' an element 
$a \in \cA$ with respect to $\phitild$ (the scalars available are 
from $\bC$, and there may be no $\lambda \in \bC$ such that 
$\phitild ( a - \lambda 1_{\cA} ) = 0$). This limits the scope 
of condition (\ref{eqn:4.81}), and makes it insufficient for 
recomputing $\phitild$ on Alg$( \cA_1 \cup \cdots \cup \cA_k )$ 
from the restrictions $\phitild \mid \cA_i$, $1 \leq i \leq k$.

For a simple concrete example showing how (\ref{eqn:4.81}) may
fail to imply infinitesimal freeness, suppose we are in the 
situation from Example \ref{ex:2.7}, with
$\cA = \bC [ \bZ_2 * \cdots * \bZ_2 ]$ and where 
$\cA_1 = \mbox{span} \{ 1_{\cA} , u_1 \}, \ldots$, 
$\cA_k = \mbox{span} \{ 1_{\cA} , u_k \}$ are the $k$ copies
of $\bC [ \bZ_2 ]$ canonically embedded into $\cA$. 
Suppose moreover that the linear functionals 
$\varphi , \phiprim : \cA \to \bC$ are such that 
$\phitild = \varphi + \ee \phiprim$ satisfies
\begin{equation}    \label{eqn:4.82}
\phitild ( 1_{\cA} ) = 1, \ \
\phitild ( u_1 ) = \cdots = \phitild ( u_k ) = \ee .
\end{equation}
Then, no matter how $\phitild$ acts on words of length 
$\geq 2$ made with $u_1, \ldots , u_k$, it will be true that 
$\cA_1, \ldots , \cA_k$ satisfy condition (\ref{eqn:4.81}) with 
respect to $\phitild$; this is due to the simple reason that the 
restrictions $\phitild \mid \cA_i$ ($1 \leq i \leq k$) are 
one-to-one. But on the other hand, Remark \ref{rem:2.2} tells 
us that if $\cA_1, \ldots , \cA_k$ are to be infinitesimally free 
in $( \cA , \varphi , \phiprim )$, then $\phitild$ is uniquely
determined by (\ref{eqn:4.82}); for example, the formulas given 
for illustration in Equations (\ref{eqn:2.21}), (\ref{eqn:2.22}) 
imply that we must have
$\phitild (u_1 u_2)  = \phitild (u_1) \phitild (u_2) 
 = \ee^2  =  0$. 
Hence any choice of $\phitild$ as in (\ref{eqn:4.82}) and with 
$\phitild (u_1 u_2) \neq 0$ provides an example for how condition 
(\ref{eqn:4.81}) does not imply infinitesimal freeness. 
\end{remark}

We conclude this section by establishing the fact about 
traciality that was announced in Remark \ref{rem:2.3}.

\begin{lemma}   \label{lemma:4.9}
Let $\cB$ be a unital subalgebra of $\cA$, and suppose that 
$\phitild \mid \cB$ is a trace. Then
\begin{equation}     \label{eqn:4.91}
\ktild_n (b_1, b_2, \ldots,  b_n ) = 
\ktild_n (b_2, b_n, \ldots , b_1), \ \ 
\forall \, n \geq 2, \ b_1, \ldots , b_n \in \cB .
\end{equation}
\end{lemma}

\begin{proof} Let $\Gamma$ be the cyclic permutation of 
$\{ 1, \ldots , n \}$ defined by $\Gamma (1) =2, \ldots ,
\Gamma (n-1) = n, \Gamma (n) = 1$. It is easy to see 
(cf. Exercise 9.41 on p. 153 of \cite{NS2006}) that $\Gamma$ 
induces an automorphism of the lattice $NC(n)$ which maps 
$\pi = \{ V_1, \ldots , V_p \} \in NC(n)$ to 
$\Gamma \cdot \pi := \{ \Gamma (V_1), \ldots , \Gamma (V_p) \}$.

Now let some $b_1, \ldots , b_n \in \cB$ be given. The right-hand
side of (\ref{eqn:4.91}) is 
$\ktild_n ( b_{\Gamma (1)}, \ldots , b_{\Gamma (n)} )$, which 
is by definition equal to 
\begin{equation}    \label{eqn:4.92}
\sum_{\pi \in NC(n)} \ \moebA ( \pi , 1_n ) \cdot
\phitild_{\pi} ( b_{\Gamma (1)}, \ldots , b_{\Gamma (n)} ).
\end{equation}
By taking into account the traciality of $\phitild$ on $\cB$
it is easily verified that 
$\phitild_{\pi} ( b_{\Gamma (1)}, \ldots , b_{\Gamma (n)} )$
= $\phitild_{\Gamma \cdot \pi} ( b_1, \ldots , b_n )$,
$\forall \, \pi \in NC(n)$. Since $\moebA ( \Gamma \cdot \pi , 1_n )$
= $\moebA ( \Gamma \cdot \pi , \Gamma \cdot 1_n ) 
= \moebA ( \pi , 1_n ), \ \ \forall \, \pi \in NC(n)$,
it becomes clear that the change of variable 
$\Gamma \cdot \pi =: \rho$ will convert the sum from (\ref{eqn:4.92})
into the one which defines $\ktild_n (b_1, \ldots , b_n)$.
\end{proof}

\begin{proposition}   \label{prop:4.10}
Let $\cA_1, \ldots , \cA_k$ be unital subalgebras of $\cA$ that 
are infinitesimally free in $( \cA , \varphi , \phiprim )$. If 
$\varphi \mid \cA_i$ and $\phiprim \mid \cA_i$ are traces for 
every $1 \leq i \leq k$, then $\varphi$ and $\phiprim$ are 
traces on Alg$( \cA_1 \cup \cdots \cup \cA_k )$.
\end{proposition}

\begin{proof} The given hypothesis and the required conclusion can
be rephrased by saying that $\phitild$ is a trace on every $\cA_i$,
and respectively that $\phitild$ is a trace on 
Alg$( \cA_1 \cup \cdots \cup \cA_k )$. Clearly, the rephrased 
conclusion will follow if we prove that 
\begin{equation}   \label{eqn:4.101}
\phitild (x_1 \cdots x_{n-1} x_n )
= \phitild (x_n x_1 \cdots x_{n-1})
\end{equation}
where $x_1 \in \cA_{i_1}, \ldots , x_n \in \cA_{i_n}$ with
$n \geq 2$ and $1 \leq i_1, \ldots , i_n \leq k$.
Let us fix such $n$, $i_1, \ldots , i_n$ and $x_1, \ldots , x_n$.
It is moreover convenient to denote
$y_1 := x_n, y_2 := x_1, \ldots , y_n := x_{n-1}$,
so that (\ref{eqn:4.101}) takes the form 
$\phitild (x_1 \cdots x_n )$ = 
$\phitild (y_1 \cdots y_n)$.

Let $\pi_o$ be the partition of $\{ 1, \ldots , n \}$ defined 
by the requirement that for $1 \leq p<q \leq n$ we have
$\bigl( \mbox{ $p,q$ in the same block of $\pi_o$ } \bigr)
\ \Leftrightarrow \ i_p = i_q$.
The hypothesis that $\cA_1, \ldots , \cA_k$ are infinitesimally 
free and Proposition \ref{prop:4.6} imply that 
\begin{equation}   \label{eqn:4.102}
\phitild (x_1 \cdots x_n ) = \sum_{ \begin{array}{c}
{\scriptstyle \pi \in NC(n) \ such } \\
{\scriptstyle that \ \pi \leq \pi_o}
\end{array} } \ \ktild_{\pi} (x_1, \ldots , x_n).
\end{equation}
(Note that $\pi_o$ may not belong to $NC(n)$, but the inequality
$\pi \leq \pi_o$ still makes sense, in reverse refinement order.)
Now, by using Lemma \ref{lemma:4.9} it is easily checked that for 
every $\pi \in NC(n)$ such that $\pi \leq \pi_o$ one has 
\begin{equation}   \label{eqn:4.103}
\ktild_{\pi} (x_1, \ldots , x_n) =
\ktild_{\Gamma \cdot \pi} (y_1, \ldots , y_n),
\end{equation}
where ``$\Gamma \cdot \pi$'' has the same significance as in the 
proof of Lemma \ref{lemma:4.9}. If we combine (\ref{eqn:4.102}) 
with (\ref{eqn:4.103}) and then make the change of variable 
$\Gamma \cdot \pi =: \rho$, we arrive to
\begin{equation}   \label{eqn:4.104}
\phitild (x_1 \cdots x_n ) = \sum_{ \begin{array}{c}
{\scriptstyle \rho \in NC(n) \ such } \\
{\scriptstyle that \ \rho \leq \Gamma \cdot \pi_o}
\end{array} } \ \ktild_{\rho} (y_1, \ldots , y_n).
\end{equation}
Finally, we invoke once more the infinitesimal freeness of 
$\cA_1, \ldots , \cA_k$ and Proposition \ref{prop:4.6}, to conclude
that the right-hand side of (\ref{eqn:4.104}) is precisely the 
moment-cumulant expansion for $\phitild (y_1 \cdots y_n )$.
\end{proof}

$\ $

$\ $

\begin{center}
{\bf\large 5. Alternating products of infinitesimally free 
random variables}
\end{center}
\setcounter{section}{5}
\setcounter{equation}{0}
\setcounter{theorem}{0}

In Proposition \ref{prop:4.6} we saw that infinitesimal freeness
can be described as a vanishing condition for mixed $\bG$-valued
cumulants. Because of this fact and because $\bG$ is commutative, 
(which makes practically all calculations with non-crossing 
cumulants go without any change from $\bC$-valued to $\bG$-valued
framework) we get a ``generic method'' for proving infinitesimal
versions of various results presented in the monograph
\cite{NS2006} -- replace $\bC$ by $\bG$ in the proof of the 
original result, then take the soul part of what comes out. Note 
that the infinitesimal results so obtained do not have $\bG$ in 
their statement, hence could also be attacked by using other 
approaches to infinitesimal freeness (in which case, however, 
proving them may be more than a straightforward routine).

In this section we show how the generic method suggested above  
works when applied to the topic of alternating products of 
infinitesimally free random variables. In particular, we will 
obtain the infinitesimal versions for two important facts 
related to this topic, that were originally found in 
\cite{NS1996} -- one of them is about compressions by free 
projections, the other concerns a method of constructing 
free families of free Poisson elements. Since the proofs of the 
$\bG$-valued formulas that we need are identical to those of
their $\bC$-valued counterparts, we will not give them here, 
but we will merely indicate where in \cite{NS2006} can the 
$\bC$-valued proofs be exactly found. The starting point is
provided by the following formulas, obtained by doing the 
$\bC$-to-$\bG$ change in Theorem 14.4 of \cite{NS2006}.

\begin{proposition}    \label{prop:5.1}
Let $( \cA , \varphi , \phiprim )$ be an incps and let 
$\cA_1, \cA_2$ be unital subalgebras of $\cA$ which are 
infinitesimally free. Consider the functional 
$\phitild = \varphi + \ee \phiprim : \cA \to \bG$ and the 
associated cumulant functionals 
$( \ktild_n : \cA^n \to \bG )_{n \geq 1}$. Recall that for every 
$n \geq 1$ and $\pi \in NC(n)$ we also have functionals 
$\phitild_{\pi}, \ktild_{\pi} : \cA^n \to \bG$, as defined in 
Notation \ref{def:3.22}.

$1^o$ For every $a_1, \ldots , a_n \in \cA_1$ and 
$b_1, \ldots , b_n \in \cA_2$ one has that 
\begin{equation}    \label{eqn:5.11}
\phitild (a_1 b_1 \cdots a_n b_n) = \sum_{\pi \in NC(n)}
\ktild_{\pi} (a_1, \ldots , a_n) \cdot 
\phitild_{Kr ( \pi ) } (b_1, \ldots , b_n).
\end{equation}

$2^o$ For every $a_1, \ldots , a_n \in \cA_1$ and 
$b_1, \ldots , b_n \in \cA_2$ one has that 
\begin{equation}    \label{eqn:5.12}
\ktild_n (a_1 b_1, \ldots , a_n b_n) = \sum_{\pi \in NC(n)}
\ktild_{\pi} (a_1, \ldots , a_n) \cdot 
\ktild_{Kr ( \pi ) } (b_1, \ldots , b_n).
\end{equation}

\hfill $\square$
\end{proposition}

We now start on the application to free compressions.

\begin{definition}          \label{def:5.2}
Let $( \cA , \varphi, \phiprim )$ be an 
incps, and let $p \in \cA$ be an idempotent element such that 
$\varphi (p) \neq 0$. We denote $\varphi (p) =: \alpha$ and 
$\phiprim (p) = \alpha '$. The {\em compression} of 
$( \cA , \varphi, \phiprim )$ by $p$ is then defined to be the
incps $( \cB , \psi , \psiprim )$ where 
\begin{equation}        \label{eqn:5.21}
\cB := p \cA p = \{ b \in \cA \mid pb = b = bp \}
\end{equation}
and where $\psi , \psiprim : \cB \to \bC$ are defined by
\begin{equation}        \label{eqn:5.22}
\psi (b) = \frac{1}{\alpha} \varphi (b), \ \ 
\psiprim (b) = \frac{1}{\alpha} \phiprim (b) - 
 \frac{\alpha '}{\alpha^2} \varphi (b), \ \ b \in \cB .
\end{equation}
\end{definition}

\begin{remark}    \label{rem:5.3}

$1^o$ In the preceding definition, note that the Grassman 
number $\widetilde{\alpha} := \alpha + \ee \alpha '$ is 
invertible in $\bG$, with inverse 
$1 / \widetilde{\alpha}$ = 
$( 1/ \alpha ) - \ee ( \alpha ' / \alpha^2 )$. As a consequence,
the two formulas given in (\ref{eqn:5.22}) are equivalent
to the fact that the consolidated functional 
$\psitild = \psi + \ee \psiprim : \cB \to \bG$ satisfies
\begin{equation}    \label{eqn:5.31}
\psitild (b) = \frac{1}{\widetilde{\alpha}} \
\phitild (b), \ \ \ \forall \, b \in \cB .
\end{equation}   

$2^o$ If in the preceding definition $( \cA , \varphi , \phiprim )$ 
is a $*$-incps and $p$ is a projection, then by using the 
relations $p = p^* = p^2$ we immediately infer that 
$0 < \alpha \leq 1$ and $\alpha ' \in \bR$. As a consequence,
$( \cB , \psi , \psiprim )$ defined there is a $*$-incps as well. 
\end{remark}

\begin{theorem}    \label{thm:5.4}
Let $( \cA , \varphi, \phiprim )$ be an incps. Let $p \in \cA$ be 
an idempotent element such that $\varphi (p) \neq 0$. Denote 
$\varphi (p) =: \alpha$, $\phiprim (p) =: \alpha '$, and consider 
the compressed incps $( \cB , \psi , \psiprim )$ from Definition 
\ref{def:5.2}. For every $n \geq 1$ let 
$\kappa_n , \kappa_n ' : \cA^n \to \bC$ and 
$\ckn , \ckn ' : \cB^n \to \bC$ be the $n$th non-crossing
cumulant and infinitesimal cumulant functional associated to 
$( \cA , \varphi, \phiprim )$ and to $( \cB , \psi , \psiprim )$,
respectively. Let $\cX$ be a subset of $\cA$ which is 
infinitesimally free from $\{ p \}$. Then we have 
\begin{equation}    \label{eqn:5.41}
\ckn (p x_1 p, \ldots , p x_n p) = 
\frac{1}{\alpha} \kappa_n ( \alpha x_1, \ldots , \alpha x_n ),
\ \ \forall \, n \geq 1, \ x_1, \ldots , x_n \in \cX
\end{equation}
and
\begin{equation}    \label{eqn:5.42}
\left\{   \begin{array}{l}
\underline{\kappa}_{1} ' (p x_1 p) = \kappa_1 ' ( x_1 ),
\ \ \forall \, x_1 \in \cX                          \\
\ckn ' (p x_1 p, \ldots , p x_n p) = 
\frac{(n-1) \alpha '}{\alpha^2} 
\kappa_n ' ( \alpha x_1, \ldots , \alpha x_n ),
\ \ \forall \, n \geq 2, \ x_1, \ldots , x_n \in \cX .
\end{array}   \right.
\end{equation}
\end{theorem}

\begin{proof} It is easily verified that Equations (\ref{eqn:5.41}) 
and (\ref{eqn:5.42}) are the body part and respectively the soul
part for the formula
\begin{equation}    \label{eqn:5.43}
\cktild (p x_1 p, \ldots , p x_n p) = 
\alphatild^{n-1} \cdot \ktild_n ( x_1, \ldots , x_n ) \in \bG,
\ \ \forall \, n \geq 1, \ x_1, \ldots , x_n \in \cX ,
\end{equation}
where the ``tilde'' notations have their usual meaning
($\cktild = \ckn + \ee \cdot \ckn '$,
$\alphatild = \alpha + \ee \cdot \alpha '$). But the latter 
formula is just the $\bG$-valued counterpart for Theorem 14.10 in
\cite{NS2006}; its proof is obtained by faithfully doing the
$\bC$-to-$\bG$ transcription of the proof of that theorem 
in \cite{NS2006}, with the minor change that the powers of 
$\alphatild$ must be kept outside the cumulant functionals 
(one cannot write 
``$\ktild_n ( \alphatild x_1, \ldots , \alphatild x_n )$'', since 
$\cA$ is only a $\bC$-algebra). Note that the argument obtained 
in this way is indeed an application of Proposition \ref{prop:5.1},
in the same way as Theorem 14.10 is an application of Theorem 14.4
in \cite{NS2006}.
\end{proof}

\begin{corollary}    \label{cor:5.5}
Let $( \cA , \varphi, \phiprim )$ be an incps. Let $p \in \cA$ be 
an idempotent element with $\varphi (p) \neq 0$, and consider the 
compressed incps $( \cB , \psi , \psiprim )$ defined as above.
Let $\cX_1, \ldots , \cX_k$ be subsets of $\cA$ such that 
$\{ p \}, \cX_1, \ldots , \cX_k$ are 
infinitesimally free in $( \cA , \varphi, \phiprim )$.
Put $\cY_i = p \cX_i p \subseteq \cB$, $1 \leq i \leq k$.
Then $\cY_1, \ldots , \cY_k$ are infinitesimally free in 
$( \cB , \psi , \psiprim )$.
\end{corollary}

\begin{proof} This is an immediate consequence of Corollary
\ref{cor:4.7}, where the needed vanishing of mixed cumulants 
follows from the explicit formulas found in Theorem \ref{thm:5.4}.
\end{proof}

We now go to the construction of families of infinitesimally free 
Poisson elements. We will use the infinitesimal (a.k.a ``type B'') 
versions of semicircular and of free Poisson elements that 
appeared in \cite{P2007} in connection to limit theorems of type B,
and are discussed in detail in Sections 4 and 5 of \cite{BS2009}.
For the present paper it is most convenient to introduce these 
elements in terms of their infinitesimal cumulants, as stated in 
Definitions \ref{def:5.6} and \ref{def:5.8} below.

\begin{definition}    \label{def:5.6}
Let $( \cA , \varphi, \phiprim )$ be a $*$-incps. A selfadjoint 
element $x \in \cA$ will be called {\em infinitesimally semicircular} 
when it satisfies 
\begin{equation}   \label{eqn:5.61}
\kappa_n (x, \ldots , x) = \kappa_n ' (x, \ldots , x) = 0,
\ \ \forall \, n \geq 3.
\end{equation}
If in addition to that we also have 
\begin{equation}   \label{eqn:5.62}
\kappa_1 (x) = 0, \ \ \kappa_2 (x,x) = 1,
\end{equation}
then we will say that $x$ is a 
{\em standard} infinitesimally semicircular element.
\end{definition}

\begin{remark}   \label{rem:5.7}

$1^o$ By using the multilinearity of $\kappa_n, \kappa_n '$ and 
Proposition \ref{prop:3.24}, it is immediately seen that if $x$ 
is infinitesimally semicircular then so is 
$\alpha ( x - \beta 1_{\cA} )$ for any $\alpha > 0$ and 
$\beta \in \bR$. Moreover, leaving aside the trivial case when 
$\kappa_2 (x,x) = 0$, one can always pick $\alpha$ and $\beta$ 
so that $\alpha ( x - \beta 1_{\cA} )$ is standard.

$2^o$ Let $x$ be standard infinitesimally semicircular in 
$( \cA , \varphi , \phiprim )$. Then all moments $\varphi (x^n)$ 
and $\phiprim (x^n)$ for $n \geq 1$ are completely determined by 
the real parameters 
\footnote{ Any two numbers $\alpha_1 ', \alpha_2 ' \in \bR$ can 
appear here. Indeed, Example \ref{ex:8.7} shows situations where 
one has $\alpha_1 ' =1, \alpha_2 ' = 0$ and respectively 
$\alpha_1 ' =0, \alpha_2 ' = 2$. One can rescale the
functionals $\phiprim$ of these two special cases to get standard
infinitesimal semicirculars $x_1, x_2$ having any pairs of parameters 
$\alpha_1 ' , 0$ and respectively $0, \alpha_2 '$; then due to 
Proposition \ref{prop:2.4} one may assume that $x_1, x_2$ are 
infinitesimally free, and form the average $(x_1 + x_2)/ \sqrt{2}$,
which is standard infinitesimally semicircular with generic 
parameters in (\ref{eqn:5.71}). }
$\alpha_1 ', \alpha_2 '$ defined by
\begin{equation}   \label{eqn:5.71}
\alpha_1 ':= \kappa_1 ' (x) = \phiprim (x), \ \ \mbox{ and } \ \ 
\alpha_2 ':= \kappa_2 ' (x,x) = \phiprim (x^2).
\end{equation}
It is in fact very easy to calculate what these moments 
are. Indeed, one can calculate the the $\bG$-valued moments 
$\phitild (x^n) = \varphi (x^n) + \ee \phiprim (x^n)$ by using the 
moment-cumulant formula (\ref{eqn:3.224}), where one takes into 
account that 
\[
\ktild_1 (x) = \ee \alpha_1 ', \
\ktild_2 (x,x) = 1 + \ee \alpha_2 ', \mbox{ and }
\ktild_n (x, \ldots , x) = 0 \mbox{ for all $n \geq 3$.}
\]
The expansion of $\phitild (x^n)$ in terms of 
$\{ \ktild_{\pi} (x, \ldots , x) \mid \pi \in NC(n) \}$ can 
get non-zero contributions only from such 
partitions $\pi$ where every block $V$ of $\pi$ has 
$|V| \leq 2$ and where there is at most one block of $\pi$ of
cardinality 1 (the latter condition coming from the fact 
that $( \ktild_1 (x) )^2 = 0$). We distinguish two cases,
depending on the parity of $n$.

{\em Case 1. $n$ is even, $n=2m$.} We get a sum extending over 
non-crossing pairings in $NC(n)$, which gives us
\[
\phitild (x^{2m}) = C_m \cdot (1+ \ee \alpha_2 ')^m  = 
C_m \cdot (1 + \ee m \alpha_2 '),
\]
or in other words 
\begin{equation}   \label{eqn:5.72}
\varphi (x^{2m}) = C_m, \ \ 
\phiprim (x^{2m}) = \alpha_2 ' \cdot (m C_m), 
\end{equation}
where $C_m$ stands for the $m$th Catalan number.

{\em Case 2. $n$ is odd, $n=2m +1$.} Here we get a sum extending
over the partitions $\pi \in NC(n)$ which have one block of 1 
element and $m$ blocks of 2 elements. There are $(2m+1) C_m$ 
such partitions; so we obtain
\[
\phitild (x^{2m+1}) = (2m+1) C_m \cdot 
\Bigl( \, ( \ee \alpha_1 ') \, (1+ \ee \alpha_2 ')^m \, \Bigr),
\]
leading to
\begin{equation}    \label{eqn:5.73}  
\varphi (x^{2m+1}) = 0, \ \ 
\phiprim (x^{2m+1}) = \alpha_1 ' \cdot \bigl( (2m+1) C_m \bigr). 
\end{equation}
\end{remark}

\begin{definition}     \label{def:5.8}
Let $( \cA , \varphi, \phiprim )$ be a $*$-incps, and let 
$\lambda , \beta ' , \gamma '$ be real parameters, where 
$\lambda > 0$. A selfadjoint element $y \in \cA$ will be called 
{\em infinitesimally free Poisson} of parameter $\lambda$ and 
\footnote{A more complete definition of these elements would also 
use a 4th parameter $r >0$, and have each of 
$\lambda, \beta ', \gamma '$ multiplied by $r^n$ in Equations 
(\ref{eqn:5.81}). For the sake of simplicity, here we have set 
this additional parameter to $r=1$.}
infinitesimal parameters $\beta ' , \gamma '$ when it has 
non-crossing cumulants given by
\begin{equation}     \label{eqn:5.81}
\left\{   \begin{array}{lclr}
\kappa_n (y, \ldots ,y)   & = &  \lambda ,             &      \\
\kappa_n ' (y, \ldots ,y) & = &  \beta ' + n \gamma ', & 
                                           \forall \, n \geq 1.
\end{array}  \right.
\end{equation}
\end{definition}

\begin{theorem}    \label{thm:5.9}
Let $( \cA , \varphi, \phiprim )$ be a $*$-incps. Let $x \in \cA$ 
be a standard infinitesimally semicircular element, and let $\cS$ 
be a subset of $\cA$ which is infinitesimally free from $\{ x \}$. 
Then for every $n \geq 1$ and $a_1, \ldots , a_n \in \cS$ we have 
\begin{equation}    \label{eqn:5.91}
\kappa_n (x a_1 x, \ldots , x a_n x) = 
\varphi (a_1 \cdots a_n)
\end{equation}
and
\begin{equation}    \label{eqn:5.92}
\kappa_n ' (x a_1 x, \ldots , x a_n x) = 
\phiprim (a_1 \cdots a_n) + 
n \, \phiprim (x^2) \cdot \varphi (a_1 \cdots a_n).
\end{equation}
\end{theorem}

\begin{proof} Equations (\ref{eqn:5.91}) and (\ref{eqn:5.92}) are 
the body part and respectively the soul part for the formula
\begin{equation}    \label{eqn:5.93}
\ktild_n (x a_1 x, \ldots , x a_n x) = 
\bigl( \ktild_2 (x,x) \bigr)^n \cdot
\phitild (a_1 \cdots a_n) \in \bG.
\end{equation}
The proof of the latter formula is obtained by doing the 
$\bC$-to-$\bG$ transcription either for the arguments used in 
Proposition 12.18 and Example 12.19 on pp. 207-208 of \cite{NS2006},
or for the arguments in Propositions 17.20 and 17.21 on pp.
283-284 of \cite{NS2006}.
\end{proof}

The ensuing construction of families of infinitesimally free 
Poisson elements is stated in the next corollary. Part $2^o$ of 
the corollary has also appeared as Corollary 36 of \cite{BS2009}.

\begin{corollary}    \label{cor:5.10}
Let $( \cA , \varphi, \phiprim )$ be a $*$-incps, and let 
$x \in \cA$ be a standard infinitesimally semicircular
element. Let $e_1, \ldots , e_k \in \cA$ be projections 
such that $e_i \perp e_j$ for $1 \leq i < j \leq k$ and such 
that $\{ e_1, \ldots , e_k \}$ is infinitesimally free from 
$\{ x \}$. Then 

$1^o$ The elements $x e_1 x, \ldots , x e_k x$ form an 
infinitesimally free family in $( \cA , \varphi , \phiprim )$.

$2^o$ For every $1 \leq i \leq k$, $x e_i x$ is infinitesimally 
free Poisson with parameter $\lambda_i$ and infinitesimal 
parameters $\beta_i ', \gamma_i '$ given by
$\lambda_i = \varphi (e_i), \ \
\beta_i ' = \phiprim (e_i), \ \ 
\gamma_i ' = \phiprim (x^2) \cdot \varphi (e_i)$. 
\end{corollary}

\begin{proof} $1^o$ This is an immediate consequence of Corollary
\ref{cor:4.7}, where the needed vanishing of mixed cumulants 
follows from the explicit formulas found in Theorem \ref{thm:5.9}.

$2^o$ By putting $a_1 = \cdots = a_n := e_i$ in (\ref{eqn:5.91}) 
and (\ref{eqn:5.92}) we see that the cumulants of $x e_i x$ have
the form required in Definition \ref{def:5.8}, with parameters 
$\lambda_i , \beta_i ', \gamma_i '$ as stated.
\end{proof}

$\ $

$\ $

\begin{center}
{\bf\large 6. Relations with the lattices \boldmath{$\ncb (n)$} }
\end{center}
\setcounter{section}{6}
\setcounter{equation}{0}
\setcounter{theorem}{0}

In this section we remember that the concept of incps has its
origins in the considerations ``of type B'' from \cite{BGN2003},
and we look at how the essence of these considerations persists 
in the framework of the present paper.

The strategy of \cite{BGN2003} was to study the type B analogue
for an operation with power series called {\em boxed convolution}
and denoted by $\freestar$. The focus on $\freestar$ was 
motivated by the fact that it provides in some sense a ``middle
ground'' between alternating products of free random variables and
the structure of intervals in the lattices $NC(n)$ (see 
discussion on pp. 2282-2283 of \cite{BGN2003}). The key point 
discovered in \cite{BGN2003} (stated in the form of the equation
$\freestar^{ \, (B)} = \freestar^{ \, (A)}_{\, \bG}$ in the 
introduction of that paper) was that boxed convolution of type B
can still be defined by the formulas from type A, provided that
one uses scalars from $\bG$.

For a detailed discussion on $\freestar$ we refer the reader to 
Lecture 17 of \cite{NS2006}. What is important for us here is that
the formula used to define $\freestar$ (cf. Equation (17.1) on 
p. 273 of \cite{NS2006}) has already made an appearance, in 
$\bG$-valued context, in Equations (\ref{eqn:5.11}),
(\ref{eqn:5.12}) of the preceding section. So then, the present
incarnation of the 
``$\, \freestar^{ \, (B)} = \freestar^{ \, (A)}_{\, \bG} \,$''
principle from \cite{BGN2003} should just amount to the following
fact: if one takes the soul parts of Equations (\ref{eqn:5.11}) 
and (\ref{eqn:5.12}), then summations over $\ncb (n)$ must arise.
This is stated precisely in Theorem \ref{thm:6.4} below, which 
is actually an easy application of the fact that the absolute 
value map $\abs : \ncb (n) \to NC(n)$ is an $(n+1)$-to-$1$ cover.

We start by introducing some notations that will be used in 
Theorem \ref{thm:6.4}, namely the type B analogues for the 
functionals $\phiA_{\pi}$ and $\kappaA_{\pi}$ from subsection 3.2.

\begin{notation}   \label{def:6.1}
Let $( \cA , \varphi , \phiprim )$ be an incps and consider the 
families of non-crossing cumulant functionals 
$( \kappa_n , \kappa_n ' )_{n \geq 1}$. For every $n \geq 1$ 
and $\tau \in \ncb (n)$ define a multilinear functional 
$\kappaB_{\tau} : \cA^n \to \bC$, as follows.

{\em Case 1.} If $\tau \in \nczb (n)$,
$\tau = \{ Z, V_1, -V_1, \ldots , V_p, - V_p \}$, then we put
\begin{equation}    \label{eqn:6.11}
\kappaB_{\tau} ( a_1, \ldots , a_n ) := {\kappa '}_{ |Z|/2 } 
\bigl( \, a_1, \ldots , a_n) \mid \abs (Z) \, \bigr) 
\cdot \prod_{j=1}^p \kappa_{ |V_j| } 
\bigl( \, (a_1, \ldots , a_n) \mid \abs (V_j) \, \bigr),
\end{equation}
for every $a_1, \ldots , a_n \in \cA$.

{\em Case 2.} If $\tau \in \ncb (n) \setminus \nczb (n)$,
$\tau = \{ V_1, -V_1, \ldots , V_p, - V_p \}$, then we put
\begin{equation}    \label{eqn:6.12}
\kappaB_{\tau} ( a_1, \ldots , a_n ) := \prod_{j=1}^p \kappa_{ |V_j| } 
\bigl( \, (a_1, \ldots , a_n) \mid \abs (V_j) \, \bigr),
\end{equation}
for $a_1, \ldots , a_n \in \cA$.
\end{notation}

\begin{notation}   \label{def:6.2}
Let $( \cA , \varphi , \phiprim )$ be an incps. Consider the 
families of multilinear functionals 
$( \varphi_n , \phiprim_n : \cA^n \to \bC )_{n \geq 1}$ defined by 
$\varphi_n = \varphi \circ \mult_n$,
$\phiprim_n = \phiprim \circ \mult_n$, where 
$\mult_n : \cA^n \to \cA$ is the multiplication map, $n \geq 1$
(same as used in Remark \ref{rem:3.14} above). 
Then for every $n \geq 1$ and every $\tau \in \ncb (n)$ we define 
a multilinear functional $\phiB_{\tau} : \cA^n \to \bC$ by the 
same recipe as in Notation \ref{def:6.1} (with discussion separated
in 2 cases), where every occurrence of $\kappa_m$ (respectively 
$\kappa_m '$) is replaced by $\varphi_m$ (respectively 
$\phiprim_m$). For example, the analogue of Case 1 is like this:
for $n \geq 1$ and for 
$\tau = \{ Z, V_1, -V_1, \ldots , V_p, - V_p \}$ in $\nczb (n)$ 
we define $\phiB_{\tau} \cA^n \to \bC$ by putting
\begin{equation}    \label{eqn:6.22}
\phiB_{\tau} ( a_1, \ldots , a_n ) := \phiprim_{ |Z|/2 } 
\bigl( \, a_1, \ldots , a_n) \mid \abs (Z) \, \bigr) 
\cdot \prod_{j=1}^p \varphi_{ |V_j| } 
\bigl( \, (a_1, \ldots , a_n) \mid \abs (V_j) \, \bigr),
\end{equation}
for $a_1, \ldots , a_n \in \cA$.
\end{notation}

\begin{remark}   \label{rem:6.3}
$1^o$ It is immediate that for 
$\tau \in \ncb (n) \setminus \nczb (n)$ one has
\begin{equation}    \label{eqn:6.31}
\kappaB_{\tau}  = \kappaA_{Abs ( \tau )}, \ \ 
\phiB_{\tau}  = \phiA_{Abs ( \tau )}.
\end{equation}

$2^o$ The functionals introduced in Notation \ref{def:6.1} 
extend both families $\kappa_n$ and $\kappa_n '$. Indeed, we
have that $\kappa_n ' = \kappaB_{1_{\pm n}}$ and that
$\kappa_n = \kappaA_{1_n} = \kappaB_{\tau}$
for every $n \geq 1$ and any $\tau \in \ncb (n)$ such that 
$\abs ( \tau ) = 1_n$ (e.g. $\tau = \{ \, \{ 1, \ldots , n \} , \
\{ -1, \ldots , -n \} \, \}$).
A similar remark holds in connection to the functionals 
$\phiB_{\tau}$ -- they extend both families $\varphi_n$ and 
$\phiprim_n$.
\end{remark}

\begin{theorem}   \label{thm:6.4}
Let $( \cA , \varphi , \phiprim )$ be an incps, and consider 
multilinear functionals on $\cA$ as in Notations \ref{def:6.1},
\ref{def:6.2}. Let $\cA_1, \cA_2$ be unital subalgebras of 
$\cA$ which are infinitesimally free. Then for every 
$a_1, \ldots , a_n \in \cA_1$ and $b_1, \ldots , b_n \in \cA_2$ 
one has
\begin{equation}   \label{eqn:6.41}
\phiprim (a_1 b_1 \cdots a_n b_n) = 
\sum_{ \sigma \in \ncb (n) } \
\kappaB_{\sigma} (a_1, \ldots , a_n) \cdot
\phiB_{Kr ( \sigma )} (b_1, \ldots , b_n) 
\end{equation}
and 
\begin{equation}   \label{eqn:6.42}
\kappa_n ' (a_1 b_1, \ldots , a_n b_n) = 
\sum_{ \sigma \in \ncb (n) } \
\kappaB_{\sigma} (a_1, \ldots , a_n) \cdot
\kappaB_{Kr ( \sigma )} (b_1, \ldots , b_n).
\end{equation}
\end{theorem}

\begin{proof} Consider the ``tilde'' notations from 
Proposition \ref{prop:5.1}. Let $\pi$ be a partition in $NC(n)$, 
and let us look at the expression
\[
\so \Bigl( \, \ktild_{\pi} (a_1, \ldots , a_n)
\, \ktild_{Kr ( \pi )} (b_1, \ldots , b_n) \, \Bigr)
\]
\[
= \so \Bigl( \, \prod_{V \in \pi} 
\ktild_{|V|} \bigl( (a_1, \ldots , a_n) \mid V \bigr) \cdot
\prod_{W \in Kr (\pi)} \ktild_{|W|} 
\bigl( (b_1, \ldots , b_n) \mid W \bigr) \, \Bigr).
\]
In view of the formula (\ref{eqn:3.211}) describing the soul part of 
a product, the latter expression is equal to a sum of $n+1$ terms,
some of them indexed by the blocks $V \in \pi$, and the others 
indexed by the blocks $W \in \Kr ( \pi )$. We leave it as a 
straightforward exercise to the reader to write these $n+1$ terms
explicitly, and verify that the natural correspondence to the $n+1$ 
partitions in $\{ \tau \in \ncb (n) \mid \abs ( \tau ) = \pi \}$
leads to the formula 
\begin{equation}   \label{eqn:6.43}
\so \Bigl( \, \ktild_{\pi} (a_1, \ldots , a_n)
\, \ktild_{Kr ( \pi )} (b_1, \ldots , b_n) \, \Bigr)
\end{equation}   
\[
= \sum_{  \begin{array}{c}
{\scriptstyle  \tau \in \ncb (n) \ such}   \\
{\scriptstyle  that \ \abs ( \tau ) = \pi } 
\end{array} } \ 
\kappaB_{\tau} (a_1, \ldots , a_n) \cdot
\phiB_{Kr ( \tau )} (b_1, \ldots , b_n).
\]
(Note: the Kreweras complement $\Kr ( \tau )$ from (\ref{eqn:6.43}) 
is taken in the lattice $\ncb (n)$; we use here the fact that 
$\abs ( \tau ) = \pi \Rightarrow \abs ( \Kr ( \tau ) ) = \Kr ( \pi )$
-- cf. Lemma 1.4 in \cite{BGN2003}.)

By summing over $\pi \in NC(n)$ on both sides of (\ref{eqn:6.43}),
we obtain that 
\[
\so \Bigl( \, \mbox{right-hand side of Equation 
(\ref{eqn:5.11})} \, \Bigr) = 
\mbox{(right-hand side of Equation (\ref{eqn:6.43}))}.
\]
Since the soul part of the left-hand side of Equation
(\ref{eqn:5.11}) is $\phiprim (a_1 b_1 \cdots a_n b_n)$, this 
proves that (\ref{eqn:6.41}) holds. The verification of 
(\ref{eqn:6.42}) is done in exactly the same way, by starting 
from Equation (\ref{eqn:5.12}) of Proposition \ref{prop:5.1}.
\end{proof}

\begin{remark}   \label{rem:6.5}
If in the preceding theorem we make $\cA_1 = \cA$ and 
$\cA_2 = \bC 1_{\cA}$, and if we take 
$b_1 = \cdots = b_n = 1_{\cA}$, then we obtain the formula
\begin{equation}   \label{eqn:6.51}
\phiprim  (a_1 \cdots a_n) = \sum_{\sigma \in \nczb (n)} 
\ \kappaB_{\sigma} (a_1, \ldots , a_n),
\ \ \forall \, a_1, \ldots , a_n \in \cA.
\end{equation}
The terms indexed by $\sigma \in \ncb (n) \setminus \nczb (n)$
have disappeared in (\ref{eqn:6.51}), due to the fact that 
$\phiprim ( 1_{\cA} ) = 0$. This formula was also noticed (via 
a direct argument from the definition of the $\bG$-valued 
functionals $\ktild_n$) in Proposition 7.4.4 of \cite{O2007}.
\end{remark}

$\ $

\begin{center}
{\bf\large 7. Dual derivation systems}
\end{center}
\setcounter{section}{7}
\setcounter{equation}{0}
\setcounter{theorem}{0}

\begin{notation}   \label{def:7.1}
Let $\cA$ be a unital algebra over $\bC$, and for every $n \geq 1$ 
let $\fM_n$ denote the vector space of multilinear functionals from 
$\cA^n$ to $\bC$. If $\pi = \{ V_1, \ldots , V_p \}$ is a partition 
in $NC(n)$ where the blocks $V_1, \ldots , V_p$ are listed in 
increasing order of their minimal elements, then we define a 
multilinear map 
\begin{equation}     \label{eqn:7.11}
J_{\pi} : \fM_{|V_1|} \times \cdots \times \fM_{|V_p|} \ni 
(f_1, \ldots , f_p) \to f \in \fM_n,
\end{equation}   
where 
\begin{equation}     \label{eqn:7.12}
f(a_1, \ldots , a_n) := \prod_{j=1}^p 
f_j \bigr( \, (a_1, \ldots , a_n) \mid V_j \, \bigr), \ \ \forall 
\, a_1, \ldots , a_n \in \cA.
\end{equation}
\end{notation}

\begin{remark}   \label{rem:7.10}

$1^o$ The formula (\ref{eqn:7.12}) from the preceding notation is 
the same as those used to define the families of functionals 
$\{ \phiA_{\pi} \mid \pi \in NC(n) \}$ and
$\{ \kappaA_{\pi} \mid \pi \in NC(n) \}$ in Remark \ref{rem:3.14}.
Hence if $( \cA , \varphi )$ is a noncommutative probability space
and if $( \kappa_n )_{n \geq 1}$ are the non-crossing cumulant 
functionals associated to $\varphi$, then for 
$\pi = \{ V_1, \ldots , V_p \} \in NC(n)$ as in Notation 
\ref{def:7.1} we get that 
\begin{equation}   \label{eqn:7.101}
J_{\pi} ( \, \kappa_{|V_1|}, \ldots , \kappa_{|V_p|} \, ) 
= \kappaA_{\pi}.
\end{equation}
Likewise, for the same $( \cA , \varphi )$ and $\pi$ we get
\begin{equation}   \label{eqn:7.102}
J_{\pi} ( \, \varphi_{|V_1|}, \ldots , \varphi_{|V_p|} \, ) 
= \phiA_{\pi},
\end{equation}
where $\varphi_m = \varphi \circ \mult_m : \cA^m \to \bC$,
$m \geq 1$ (same as in Remark \ref{rem:3.14}).

$2^o$ Let $\pi = \{ V_1, \ldots , V_p \} \in NC(n)$ be as in 
Notation \ref{def:7.1}, and let $1 \leq j \leq p$ be such that 
$V_j$ is an interval-block of $\pi$. Denote $|V_j| =: m$ and 
let $\piciup \in NC(n-m)$ be the partition obtained by 
removing the block $V_j$ out of $\pi$ and by redenoting the 
elements of $\{ 1, \ldots , n \} \setminus V_j$ as 
$1, \ldots , n-m$, in increasing order. On the other hand, 
let us denote by $\gamma \in NC(n)$ the partition of 
$\{ 1, \ldots , n \}$ into the two blocks $V_j$ and 
$\{ 1, \ldots , n \} \setminus V_j$. It is then immediate that 
for every $f_1 \in \fM_{|V_1|}, \ldots , f_p \in \fM_{|V_p|}$
we can write
\begin{equation}    \label{eqn:7.103}
J_{\pi} (f_1, \ldots , f_p) = J_{\gamma} (g, f_j)  
\mbox{$\ \ $ where } g := 
   J_{\piciup} (f_1, \ldots  f_{j-1}, f_{j+1}, \ldots , f_p).
\end{equation} 
Due to this observation and to the fact that every non-crossing
partitions has interval-blocks, considerations about the 
multilinear functions $J_{\pi}$ from Notation \ref{def:7.1} can 
sometimes be reduced (via an induction argument on $| \pi |$) to
discussing the case when $| \pi | = 2$.
\end{remark}

\begin{definition}   \label{def:7.2}
Let $\cA$ be a unital algebra over $\bC$ and let the spaces 
$( \fM_n )_{n \geq 1}$ and the multilinear functions 
$\{ J_{\pi} \mid \pi \in \cup_{n=1}^{\infty} NC(n) \}$ be as in 
Notation \ref{def:7.1}. We will call {\em dual derivation system} 
a family of linear maps $( d_n : \fD_n \to \fM_n )_{n \geq 1}$ 
where, for every $n \geq 1$, $\fD_n$ is a a linear subspace of 
$\fM_n$, and where the following two conditions are satisfied.

(i) Let $\pi = \{ V_1, \ldots , V_p \} \in NC(n)$ be as in 
Notation \ref{def:7.1}. Then for every 
$f_1 \in \fD_{|V_1|}, \ldots , f_p \in \fD_{|V_p|}$ one has 
that $J_{\pi} ( f_1, \ldots , f_p) \in \fD_n$ and that 
\begin{equation}   \label{eqn:7.21}
d_n \bigl( \, J_{\pi} ( f_1, \ldots , f_p) \, \bigr) = 
\sum_{j=1}^p   J_{\pi} ( f_1, \ldots , f_{j-1}, d_{|V_j|} (f_j), 
                              f_{j+1}, \ldots , f_p).
\end{equation}

(ii) For every $f \in \fD_1$ and every $n \geq 1$ one has that 
$f \circ \mult_n \in \fD_n$ and that
\begin{equation}   \label{eqn:7.22}
d_n \bigl( f \circ \mult_n ) = ( d_1 \,  f ) \circ \mult_n,
\end{equation}
where $\mult_n : \cA^n \to \cA$ is the multiplication map.
\end{definition}

\begin{remark}   \label{rem:7.3}

$1^o$ When verifying condition (i) in Definition \ref{def:7.2}, 
it suffices to check the particular case when $| \pi | = 2$. 
Indeed, the general case of Equation (\ref{eqn:7.21}) can then 
be obtained by induction on $| \pi |$, where one invokes the
argument from (\ref{eqn:7.103}).

$2^o$ In the setting of Definition \ref{def:7.2}, let us use 
the notation $f \times g$ for the functional obtained by 
``concatenating'' $f \in \fM_m$ and $g \in \fM_n$. So 
$f \times g \in \fM_{m+n}$ acts simply by
\[
( f \times g ) (a_1, \ldots , a_m, b_1, \ldots , b_n) =
f( a_1, \ldots , a_m) g( b_1, \ldots , b_n), \ \ \
\forall \, a_1, \ldots , a_m, b_1, \ldots , b_n \in \cA .
\]
Clearly one can write $f \times g = J_{\gamma} (f,g)$ where 
$\gamma \in NC(m+n)$ is the partition with two blocks 
$\{ 1, \ldots , m \}$ and $\{ m+1, \ldots , m+n \}$. By using 
Equation (\ref{eqn:7.21}) we thus obtain that 
\begin{equation}   \label{eqn:7.32}
d_{m+n} ( f \times g) = \bigl( d_m (f) \times g \bigr) + 
\bigl( f \times d_n (g) \bigr) ,  
\ \ \ \forall \, m,n \geq 1, \ f \in \fM_m, \, g \in \fM_n .
\end{equation}
So a dual derivation system gives in particular a derivation 
on the algebra structure defined by using concatenation on 
$\oplus_{n=1}^{\infty} \fM_n$. Note however that Equation 
(\ref{eqn:7.32}) alone is not sufficient to ensure condition (i)
from Definition \ref{def:7.2} (since it cannot control $J_{\pi}$ 
for partitions such as 
$\pi = \{ \, \{ 1,3 \}, \, \{ 2 \} \ \} \in NC(3)$).
\end{remark}

\begin{proposition}    \label{prop:7.4}
Let $\cA$ be a unital algebra over $\bC$ and let 
$( d_n : \fD_n \to \fM_n )_{n \geq 1}$ be a dual derivation
system on $\cA$. Let $\varphi$ be a linear functional in 
$\fD_1$, and denote $d_1 ( \varphi ) =: \phiprim$. Consider the 
incps $( \cA , \varphi , \phiprim )$, and let 
$( \kappa_n , \kappa_n ')_{n \geq 1}$ be the non-crossing 
cumulant and infinitesimal cumulant functionals associated to 
this incps. Then for every $n \geq 1$ we have that 
\begin{equation}   \label{eqn:7.61}
\kappa_n \in \fD_n \mbox{ and } d_n ( \kappa_n ) = \kappa_n '.
\end{equation}
\end{proposition}

\begin{proof} Denote as usual $\varphi_n := \varphi \circ \mult_n$,
$\phiprim_n := \phiprim \circ \mult_n$, $n \geq 1$. Since 
$\varphi \in \fD_1$, condition (ii) from Definition \ref{def:7.2} 
implies that $\varphi_n \in \fD_n$ and 
$d_n ( \varphi_n ) = \phiprim_n$ for every $n \geq 1$.

Now let $\pi = \{ V_1, \ldots , V_p \}$ be a partition in $NC(n)$,
with $V_1, \ldots , V_p$ written in increasing order of their 
minimal elements. By using Equation (\ref{eqn:7.102}) from 
Remark \ref{rem:7.10} and condition (i) in Definition \ref{def:7.2}
we find that 
\begin{equation}   \label{eqn:7.42}
d_n ( \phiA_{\pi} ) = 
\sum_{j=1}^p   J_{\pi} \bigl( \, \varphi_{|V_1|}, \ldots , 
\varphi_{|V_{j-1}|}, \phiprim_{|V_j|}, \varphi_{|V_{j+1}|}, 
\ldots , \varphi_{|V_p|} \, \bigr)
\end{equation}
(where the latter formula incorporates the fact that 
$d_{|V_j|} ( \varphi_{|V_j|} ) = \phiprim_{|V_j|}$).

We next consider the formula (\ref{eqn:3.143}) which expresses
a cumulant functional $\kappa_n$ in terms of the functionals 
$\{ \phiA_{\pi} \mid \pi \in NC(n) \}$. From this formula it 
follows that $\kappa_n \in \fD_n$ and that 
\begin{equation}    \label{eqn:7.43}
d_n ( \kappa_n ) = 
 \sum_{ \begin{array}{c}
{\scriptstyle \pi \in NC(n),} \\
{\scriptstyle \pi = \{ V_1, \ldots , V_p \} }
\end{array}  } \ \moeb ( \pi , 1_n ) 
\Bigl( \ \sum_{j=1}^p   J_{\pi} \bigl( \, \varphi_{|V_1|}, \ldots , 
\varphi_{|V_{j-1}|}, \phiprim_{|V_j|}, \varphi_{|V_{j+1}|}, 
\ldots , \varphi_{|V_p|} \, \bigr) \ \Bigr) .
\end{equation}
It is immediate that on the right-hand side of (\ref{eqn:7.43}) 
we have obtained precisely the sum over 
$\{ ( \pi , V ) \mid \pi \in NC(n)$, $V$ block of $\pi \}$ 
which was used to introduce $\kappa_n '$ in Definition 
\ref{def:4.2}.
\end{proof}

\begin{proposition}  \label{prop:7.5} 
Let $( \cA , \varphi , \phiprim )$ be an incps, and consider 
the multilinear functionals $\phiA_{\pi}$ ($\pi \in NC(n)$, 
$n \geq 1$) which were introduced in Remark \ref{rem:3.14}.
Suppose that for every $n \geq 1$ the set
$\{ \phiA_{\pi} \mid \pi \in NC(n) \}$ is linearly independent 
in $\fM_n$; let $\fD_n$ denote its span, and let
$d_n : \fD_n \to \fM_n$ be the linear map defined by the 
requirement that 
\begin{equation}   \label{eqn:7.51}
d_n ( \, \phiA_{\pi} \, ) = \sum_{ \begin{array}{c}
{\scriptstyle \tau \in \nczb (n) \ such } \\
{\scriptstyle that \ Abs ( \tau ) = \pi}
\end{array}  } \  \phiB_{\tau}, \ \ \ \forall \, \pi \in NC(n),
\end{equation}
with $\phiB_{\tau}$ as in Notation \ref{def:6.2}.
Then $(d_n)_{n \geq 1}$ is a dual derivation system, and
$d_1 ( \varphi ) = \phiprim$.
\end{proposition}

\begin{proof} 
It is obvious that the unique partition $\tau \in \nczb (n)$ 
such that $\abs ( \tau ) = 1_n$ is $\tau = 1_{\pm n}$. Thus if 
we put $\pi = 1_n$ in Equation (\ref{eqn:7.51}) we obtain that 
$d_n ( \phiA_{1_n} ) = \phiB_{1_{\pm n}}$; in other words, 
this means that 
\begin{equation}   \label{eqn:7.52}
d_n ( \varphi \circ \mult_n ) = \phiprim \circ \mult_n, \ \ 
\forall \, n \geq 1.
\end{equation}
The particular case $n=1$ of (\ref{eqn:7.52}) gives us that 
$d_1 ( \varphi ) = \phiprim$. Moreover, it becomes clear that 
\[
d_n ( f \circ \mult_n ) = (d_1 f) \circ \mult_n , \ \ 
\forall \, n \geq 1 \mbox{ and } f \in \bC \varphi ;
\]
since in this proposition we have $\fD_1 = \bC \varphi$, 
we thus see that condition (ii) from Definition \ref{def:7.2} is 
verified.

The rest of the proof is devoted to verifying (i) from Definition 
\ref{def:7.2}. We fix a partition
$\pi = \{ V_1, \ldots , V_p \} \in NC(n)$ for which we will prove
that Equation (\ref{eqn:7.21}) holds. Both sides of (\ref{eqn:7.21})
behave multilinearly in the arguments 
$f_1 \in \fD_{|V_1|}, \ldots, f_p \in \fD_{|V_p|}$; hence, due to
how $\fD_{|V_1|}, \ldots, \fD_{|V_p|}$ are defined, it suffices to
prove the following statement: for every 
$\pi_1 \in NC( |V_1| ), \ldots , \pi_p \in NC( |V_p| )$  we have
that $J_{\pi} ( \phiA_{\pi_1}, \ldots , \phiA_{\pi_p} ) \in \fD_n$
and that 
\begin{equation}   \label{eqn:7.53}
d_n \bigl( \, J_{\pi} 
( \phiA_{\pi_1}, \ldots , \phiA_{\pi_p} ) \, \bigr) = 
\end{equation}
\[
\sum_{j=1}^p  J_{\pi} ( \phiA_{\pi_1}, \ldots , \phiA_{\pi_{j-1}}, 
d_{|V_j|} ( \phiA_{\pi_j} ), \phiA_{\pi_{j+1}}, \ldots ,
\phiA_{\pi_p} ).
\]
In what follows we fix some partitions 
$\pi_1 \in NC( |V_1| ), \ldots , \pi_p \in NC( |V_p| )$, for which 
we will prove that this statement holds.

Observe that, in view of how the maps $d_{ |V_j| }$ are defined, on 
the right-hand side of (\ref{eqn:7.53}) we have
\[
\sum_{j=1}^p  \, \sum_{ \begin{array}{c} 
{\scriptstyle \tau \in \nczb (n) \ such } \\
{\scriptstyle that \ Abs ( \tau ) = \pi_j}
\end{array}  } \  
J_{\pi} ( \phiA_{\pi_1}, \ldots , \phiA_{\pi_{j-1}}, 
\phiB_{\tau}, \phiA_{\pi_{j+1}}, \ldots , \phiA_{\pi_p} ).
\]
But let us recall from Remark \ref{rem:3.5} that the partitions 
in $\{ \tau \in \nczb (n) \mid \abs ( \tau ) = \pi_j \}$ are 
indexed by the set of blocks of $\pi_j$. More precisely, for
every $1 \leq j \leq p$ and $V \in \pi_j$ let us denote by 
$\tau (j,V)$ the unique partition in $\nczb (n)$ such that 
$\abs ( \tau ) = \pi_j$ and such that the zero-block $Z$ of 
$\tau$ has $\abs (Z) = V$; then the double sum written above 
for the right-hand side of Equation (\ref{eqn:7.53}) becomes
\begin{equation}   \label{eqn:7.54}
\sum_{j=1}^p  \, \sum_{ V \in \pi_j } \
J_{\pi} ( \phiA_{\pi_1}, \ldots , \phiA_{\pi_{j-1}}, 
\phiB_{\tau (j,V)}, \phiA_{\pi_{j+1}}, \ldots , \phiA_{\pi_p} ).
\end{equation}

Now to the left-hand side of (\ref{eqn:7.53}). For every 
$1 \leq j \leq p$ let $\widehat{\pi}_j$ be the partition
of $V_j$ obtained by transporting the blocks of $\pi_j$ via the 
unique order preserving bijection from $\{ 1, \ldots , |V_j| \}$ 
onto $V_j$. Then $\widehat{\pi}_1, \ldots , \widehat{\pi}_p$ form
together a partition $\rho \in NC(n)$ which refines $\pi$, and it 
is immediate that 
$J_{\pi} ( \phiA_{\pi_1}, \ldots , \phiA_{\pi_p} ) = \phiA_{\rho}$.
In particular this shows of course that 
$J_{\pi} ( \phiA_{\pi_1}, \ldots , \phiA_{\pi_p} ) \in \fD_n$.
Moreover, by using how $d_n ( \phiA_{\rho} )$ is defined, we obtain 
that the left-hand side of (\ref{eqn:7.53}) is equal to
$\sum_{W \in \rho} \phiB_{\sigma (W)}$,
where for every $W \in \rho$ we denote by $\sigma (W)$ the 
unique partition in $\nczb (n)$ such that 
$\abs ( \sigma (W) ) = \rho$ and such that the zero-block $Z$ of 
$\sigma (W)$ has $\abs (Z) = W$.

Finally, we observe that the set of blocks of $\rho$ is the 
disjoint union of the sets of blocks of the partitions 
$\widehat{\pi}_1, \ldots , \widehat{\pi}_p$, 
and is hence in natural bijection with 
$\{ (j,V) \mid 1 \leq j \leq p$ and $V \in \pi_j \}$. We leave it 
as a straightforward (though somewhat notationally tedious) 
exercise to the reader to verify that when $W \in \rho$ corresponds 
to $(j,V)$ via this bijection, then the term indexed by $(j,V)$ in 
(\ref{eqn:7.54}) is precisely equal to $\phiB_{\sigma (W)}$. Hence 
the double sum from (\ref{eqn:7.54}) is identified term by term to
$\sum_{W \in \rho} \phiB_{\sigma (W)}$
via the bijection $W \leftrightarrow (j,V)$, and 
the required formula (\ref{eqn:7.53}) follows.
\end{proof}

\begin{remark}   \label{rem:7.6}
The linear independence hypothesis in Proposition \ref{prop:7.5} 
is necessary, otherwise we need some relations to be satisfied by 
$\varphi$ and $\phiprim$. Indeed, suppose for example that the 
set $\{ \phiA_{\pi} \mid \pi \in NC(2) \}$ is linearly dependent 
in $\fM_2$. It is immediately verified that this is equivalent to
the fact that $\varphi$ is a character of $\cA$
($\varphi (ab) = \varphi (a) \varphi (b)$, 
$\forall \, a,b \in \cA$). Hence $\kappa_2 = 0$, so if Proposition 
\ref{prop:7.5} is to work then we should have 
$\kappa_2 ' = d_2 ( \kappa_2 ) = 0$ as well, implying that 
$\phiprim$ satisfies the condition
$\phiprim (ab) = \varphi (a) \phiprim (b) + 
\phiprim (a) \varphi (b), \ \ \forall \, a,b \in \cA$.
\end{remark}

$\ $

\begin{center}
{\bf\large 8. Soul companions for a given \boldmath{$\varphi$} }
\end{center}
\setcounter{section}{8}
\setcounter{equation}{0}
\setcounter{theorem}{0}

In this section we elaborate on the facts announced in the 
subsection 1.3 of the introduction. We start by recording some 
basic properties of the set of functionals $\phiprim$ which can 
appear as soul-companions for $\varphi$, when $( \cA , \varphi )$ 
and $\cA_1, \ldots , \cA_k$ are given.

\begin{proposition}  \label{prop:8.1}
Let $( \cA , \varphi )$ be a noncommutative probability space 
and let $\cA_1, \ldots , \cA_k$ be unital subalgebras of $\cA$ 
which are freely independent in $( \cA , \varphi )$.

$1^o$ The set of linear functionals 
\begin{equation}  \label{eqn:8.11}
\cF ' := \Bigl\{ \phiprim : \cA \to \bC 
\begin{array}{cl}
\vline & \phiprim \mbox{ linear, $\phiprim ( 1_{\cA} ) = 0$,
         and } \cA_1, \ldots , \cA_k                           \\
\vline & \mbox{are infinitesimally free 
               in $( \cA , \varphi , \phiprim )$ }
\end{array} \Bigr\}
\end{equation}
is a linear subspace of the dual of $\cA$.

$2^o$ Suppose that Alg$( \cA_1 \cup \cdots \cup \cA_k )$ = $\cA$,
and consider the linear map 
\begin{equation}  \label{eqn:8.12}
\cF ' \ni \phiprim \mapsto 
( \phiprim \mid \cA_1, \ldots , \phiprim \mid \cA_k ) 
\in \cF_1 ' \times \cdots \times \cF_k ' ,
\end{equation}
where $\cF '$ is as in (\ref{eqn:8.11}) and where for 
$1 \leq i \leq k$ we denote 
$\cF_i ' = \{ \phiprim : \cA_i \to \bC \mid \phiprim$ linear, 
$\phiprim ( 1_{\cA} ) = 0 \}$. The map from (\ref{eqn:8.12}) is 
one-to-one.
\end{proposition}

\begin{proof} $1^o$ This is immediate from Definition \ref{def:1.1},
and specifically from the fact that $\phiprim$ makes a linear
appearance on the right-hand side of Equation (\ref{eqn:1.15}).

$2^o$ Let $\phiprim \in \cF '$ be such that 
$\phiprim \mid \cA_i = 0$, $\forall \, 1 \leq i \leq k$. Then 
from Equation (\ref{eqn:1.15}) it is immediate that 
$\phiprim (a_1 \cdots a_n) = 0$ for all choices of 
$a_1, \ldots , a_n \in \cA_1 \cup \cdots \cup \cA_k$. The linear 
span of the products $a_1 \cdots a_n$ formed in this way is the 
algebra generated by $\cA_1 \cup \cdots \cup \cA_k$, hence is all 
of $\cA$, and the conclusion that $\phiprim = 0$ follows.
\end{proof}

\begin{remark}   \label{rem:8.2}
In the framework of Proposition \ref{prop:8.1}, the linear map 
(\ref{eqn:8.12}) may not be surjective. For an example, consider 
the full Fock space over $\bC^2$,
\[
\cT = \bC \Omega \oplus \bC^2 \oplus ( \bC^2 \otimes \bC^2 ) 
\oplus \cdots \oplus ( \bC^2 )^{\otimes n} \oplus \cdots \ ,
\]
and let $L_1, L_2 \in B( \cT )$ be the left-creation operators
associated to the two vectors in the canonical orthonormal basis 
of $\bC^2$. Then $L_1, L_2$ are isometries with mutually orthogonal 
ranges; this is recorded in algebraic form by the relations
\[
L_1^* L_1 = L_2^* L_2 = 1 \mbox{ (identity operator on $\cT$), } 
\ \ L_1^* L_2 = 0.
\]
For $i = 1,2$ let $\cA_i$ denote the unital $*$-subalgebra of 
$B( \cT )$ generated by $L_i$, and let 
$\cA = \mbox{Alg} ( \cA_1 \cup \cA_2 )$, the unital $*$-algebra 
generated by $L_1$ and $L_2$ together. It is well-known (see 
e.g. Lecture 7 of \cite{NS2006}) that $\cA_1$ and $\cA_2$ are 
free in $( \cA , \varphi )$ where $\varphi$ is the vacuum-state 
on $\cA$. Let $\phiprim_2 : \cA_2 \to \bC$ be any linear 
functional such that $\phiprim_2 ( 1_{\cA} ) = 0$ and 
$\phiprim_2 (L_2) = 1$. Then there exists no linear functional
$\phiprim : \cA \to \bC$ such that 
$\phiprim \mid \cA_2 = \phiprim_2$ and such that $\cA_1, \cA_2$
are infinitesimally free in $( \cA , \varphi , \phiprim )$.
Indeed, if such $\phiprim$ would exist then from Equation 
(\ref{eqn:2.23}) of Remark \ref{rem:2.2} it would follow that 
\[
\phiprim ( L_1^* L_2 L_1) = \varphi (L_1^* L_1) \phiprim (L_2) 
+ \phiprim (L_1^* L_1) \varphi (L_2) 
= 1 \cdot 1 + 0 \cdot 0 = 1,
\]
which is not possible, since $L_1^* L_2 L_1 = 0$.
\end{remark}

\begin{remark}  \label{rem:8.3}
The example from the above remark shows that we can't always 
extend a given system of functionals $\phiprim_i$ in order to get 
a soul companion $\phiprim$ for $\varphi$. But Proposition 
\ref{prop:2.4} gives us an important case when we are sure 
this is possible, namely the one when $( \cA , \varphi )$ is 
the free product 
$( \cA_1 , \varphi_1 )* \cdots *( \cA_k , \varphi_k )$.

In the remaining part of this section we will look at the two 
recipes for obtaining a soul companion that were stated in 
Corollary \ref{cor:1.4} and Proposition \ref{prop:1.5}.
For the first of them, we start by verifying that a derivation on 
$\cA$ does indeed define a dual derivation system as indicated in 
Equation (\ref{eqn:1.31}).
\end{remark}

\begin{proposition}    \label{prop:8.4}
Let $\cA$ be a unital algebra over $\bC$ and let 
$D : \cA \to \cA$ be a derivation. For every $n \geq 1$ let 
$\fM_n$ denote the space of multilinear functionals from $\cA^n$ 
to $\bC$, and define $d_n : \fM_n \to \fM_n$ by putting
\begin{equation}  \label{eqn:8.41}
( d_n f ) (a_1, \ldots , a_n) := \sum_{m=1}^n
f \Bigl( a_1, \ldots , a_{m-1}, D(a_m), a_{m+1}, \ldots , 
a_n \Bigr), 
\end{equation}
for $f \in \fM_n$ and  $a_1, \ldots , a_n \in \cA$.
Then $(d_n)_{n \geq 1}$ is a dual derivation system on $\cA$.
\end{proposition}

\begin{proof} We first do the immediate verification of condition 
(ii) from Definition \ref{def:7.2}. Let $f$ be a 
functional in $\fM_1$, let $n$ be a positive integer, and 
denote $g = f \circ \mult_n \in \fM_n$. Then for every 
$a_1, \ldots , a_n \in \cA$ we have
\[
(d_n g) (a_1, \ldots , a_n) 
= \sum_{m=1}^n f \Bigl( a_1 \cdots a_{m-1} \cdot 
D(a_m) \cdot a_{m+1} \cdots a_n \Bigr)                 
= f \Bigl( \, D( a_1 \cdots a_n ) \, \Bigr)           
\]
(where at the first equality sign we used 
the definitions of $d_n$ and of $g$, and at the second 
equality sign we used the derivation property of $D$).
Since $d_1 f$ is just $f \circ D$, it is clear that we have 
obtained $d_n g = (d_1 f) \circ M_n$, as required.

For the remaining part of the proof we fix 
$\pi = \{ V_1, \ldots , V_p \} \in NC(n)$ and 
$f_1 \in \fM_{|V_1|}, \ldots , f_p \in \fM_{|V_p|}$ as in (i)
of Definition \ref{def:7.2}, and we verify that the formula 
(\ref{eqn:7.21}) holds. Denote 
$f := J_{\pi} (f_1, \ldots , f_p) \in \fM_n$. In the summation 
which defines $d_n f$ in Equation (\ref{eqn:8.41}) we group 
the terms by writing
\begin{equation}    \label{eqn:8.42}
\sum_{j=1}^p \Bigl( \, \sum_{m \in V_j}
f \Bigl( a_1, \ldots , a_{m-1}, D(a_m), a_{m+1}, \ldots , 
a_n \Bigr) \, \Bigr).
\end{equation}
It will clearly suffice to prove that, for every
$1 \leq j \leq p$, the term indexed by $j$ in the sum 
(\ref{eqn:8.42}) is equal to the term indexed by $j$ on the 
right-hand side of (\ref{eqn:7.21}).

So then let us also fix a $j$, $1 \leq j \leq p$. 
We write explicitly the block $V_j$ of $\pi$ as 
$\{ v_1, \ldots , v_s \}$ with $v_1 < \cdots < v_s$. From 
the definition of $f$ as $J_{\pi} (f_1, \ldots , f_p)$ it 
is then immediate that for $m = v_r \in V_j$ we have
\begin{equation}    \label{eqn:8.43}
f \Bigl( a_1, \ldots , a_{m-1}, D(a_m), a_{m+1}, \ldots , 
a_n \Bigr) \, \Bigr) =
\end{equation}
\[
= \Bigl( \ \prod_{  \begin{array}{c}
{\scriptstyle 1 \leq i \leq p,}  \\
{\scriptstyle i \neq j}
\end{array}  } 
\ f_i \bigl( \, (a_1, \ldots ,a_n) \mid V_i \, \bigr) \ \Bigr) 
\cdot
\ f_j \bigl( \, a_{v_1}, \ldots , a_{v_{r-1}}, D( a_{v_r} ),
a_{v_{r+1}}, \ldots , a_{v_s} \, \bigr) .
\]
When summing over $1 \leq r \leq s$ in (\ref{eqn:8.43}), the sum 
on the right-hand side only affects the last factor of the 
product, which gets summed to 
$( d_s f_j) \bigl( \, a_{v_1}, \ldots , a_{v_s} \, \bigr)$.
The result of this summation is hence that 
\[
\sum_{m \in V_j} f \Bigl( a_1, \ldots , a_{m-1}, D(a_m), 
a_{m+1}, \ldots , a_n \Bigr) \, \Bigr)
= J_{\pi} \bigl( f_1, \ldots , f_{j-1}, d_{|V_j|} (f_j), 
f_{j+1}, \ldots , f_p \bigr),
\]
as required.
\end{proof}

\begin{corollary}    \label{cor:8.5}
Let $( \cA ,\varphi )$ be a noncommutative probability space, and
let $D : \cA \to \cA$ be a derivation. Define 
$\phiprim :=\varphi \circ D$. Let the non-crossing and the 
infinitesimal non-crossing cumulant functionals associated to 
$( \cA , \varphi , \phiprim )$ be denoted by $\kk_n$ and 
respectively by $\kk_n '$, in the usual way. Then for every 
$n \geq 1$ and $a_1, \ldots , a_n \in \cA$ one has
\begin{equation}    \label{eqn:8.51}
\kappa_n '  (a_1, \ldots , a_n) = \sum_{m=1}^n 
\kappa_n \Bigl( a_1, \ldots , a_{m-1}, D(a_m), a_{m+1}, \ldots ,
a_n \Bigr) .
\end{equation}
\end{corollary}

\begin{proof} This follows from Proposition \ref{prop:7.4},
where we use the specific dual derivation system put into evidence 
in Proposition \ref{prop:8.4}.
\end{proof}

\begin{corollary}    \label{cor:8.6}
In the notations of Corollary \ref{cor:8.5}, let
$\cA_1, \ldots , \cA_k$ be unital subalgebras of $\cA$ which are 
freely independent with respect to $\varphi$, and such that
$D( \cA_i ) \subseteq \cA_i$ for $1 \leq i \leq k$.
Then $\cA_1, \ldots , \cA_k$ are infinitesimally free in 
$( \cA , \varphi , \phiprim )$.
\end{corollary}

\begin{proof}
We verify that condition (2) from Theorem \ref{thm:1.2} is 
satisfied. The vanishing of mixed cumulants $\kappa_n$ follows 
from the hypothesis that $\cA_1, \ldots , \cA_k$ are free in 
$( \cA , \varphi )$. But then the specific formula obtained for 
the infinitesimal cumulants $\kappa_n '$ in Corollary \ref{cor:8.5},
together with the hypothesis that $\cA_1, \ldots , \cA_k$ are 
invariant under $D$, imply that the mixed infinitesimal 
cumulants $\kappa_n '$ vanish as well.
\end{proof}

\begin{example}    \label{ex:8.7}
Consider the situation where $\cA$ is the algebra 
$\bC \langle X_1, \ldots , X_k \rangle$ of noncommutative 
polynomials in $k$ indeterminates. We will view $\cA$ as a 
$*$-algebra, with $*$-operation uniquely determined by the 
requirement that each of $X_1, \ldots , X_k$ is selfadjoint. 
Consider the unital $*$-subalgebras 
$\cA_1, \ldots , \cA_k \subseteq \cA$ where 
$\cA_i = \mbox{span} \{ X_i^n \mid n \geq 0 \}$,
$1 \leq i \leq k$. We will look at two natural derivations on 
$\cA$ that leave $\cA_1 , \ldots , \cA_k$ invariant, and 
we will examine some examples of infinitesimal freeness 
given by these derivations.

\vspace{6pt}

(a) Let $D : \cA \to \cA$ be the linear operator defined by 
putting $D(1) = 0$, $D(X_i) = 1$ $\forall \, 1 \leq i \leq k$, and 
\begin{equation}    \label{eqn:8.71}
D( X_{i_1} \cdots X_{i_n} ) = \sum_{m=1}^n
X_{i_1} \cdots X_{i_{m-1}} 
X_{i_{m+1}} \cdots X_{i_n}, \ \ 
\forall \, n \geq 2, \ \forall \, 1 \leq i_1, \ldots , i_n \leq k.
\end{equation}
It is immediate that $D$ is a derivation on $\cA$, which is 
selfadjoint (in the sense that $D(P^*) = D(P)^*$, 
$\forall \, P \in \cA$). For every $1 \leq i \leq k$ we have that 
$D( \cA_i ) \subseteq \cA_i$ and that $D$ acts on $\cA_i$ as the 
usual derivative (in the sense that $D( \, P(X_i) \, )$ = 
$P'(X_i)$, $\forall \, P \in \bC [X]$).

Now let $\mu : \cA \to \bC$ be a positive definite functional 
with $\mu (1) = 1$ and such that $\cA_1, \ldots , \cA_k$ are free
in $( \cA , \mu )$. Then Corollary \ref{cor:8.6} implies that 
$\cA_1, \ldots , \cA_k$ are infinitesimally free in the  
$*$-incps $( \cA , \mu , \mu ' )$, where $\mu ' := \mu \circ D$.

Note that in this special example we actually have 
\begin{equation}   \label{eqn:8.72}
\kappa_n ' (X_{i_1}, \ldots , X_{i_n}) = 0, \ \ \forall \,
n \geq 2, \ \forall \, 1 \leq i_1, \ldots ,i_n \leq k;
\end{equation}
this is an immediate consequence of the the formula (\ref{eqn:8.51}),
combined with the fact that a non-crossing cumulant vanishes when 
one of its arguments is a scalar. 

Equation (\ref{eqn:8.72}) gives in particular that
\[
\kappa_n ' (X_i, \ldots , X_i) = 0, \ \ \forall \,
n \geq 2 \mbox{ and } 1 \leq i \leq k.
\]
So if $\mu$ is defined such that every $X_i$ has a standard 
semicircular distribution in $( \cA , \mu )$, then every $X_i$ will
become a standard infinitesimal semicircular element in 
$( \cA , \mu , \mu ' )$, in the sense of Remark \ref{rem:5.7}, 
and where in Equation (\ref{eqn:5.71}) we take 
$\alpha_1 ' = 1$, $\alpha_2 ' = 0$.

\vspace{6pt}

(b) Let $\Dno : \cA \to \cA$ be the linear operator defined by 
putting $\Dno (1) = 0$ and 
\begin{equation}    \label{eqn:8.73}
\Dno ( X_{i_1} \cdots X_{i_n} ) = n \, X_{i_1} \cdots X_{i_n}, \ \ 
\forall \, n \geq 1, \ \forall \, 1 \leq i_1, \ldots , i_n \leq k.
\end{equation}
Then $\Dno$ is a selfadjoint derivation, sometimes called ``the 
number operator'' on $\cA$. It is clear that $\Dno$ leaves every
$\cA_i$ invariant, $1 \leq i \leq k$. Hence if 
$\mu : \cA \to \bC$ is as in part (a) above (such that 
$\cA_1, \ldots , \cA_k$ are free in $( \cA , \mu )$), then 
Corollary \ref{cor:8.6} implies that 
$\cA_1, \ldots , \cA_k$ are infinitesimally free in the  
$*$-incps $( \cA , \mu , \muno ' )$, where 
$\muno ' := \mu \circ \Dno$.

Since $\Dno (X_i) = X_i$ for $1 \leq i \leq k$, the formula 
(\ref{eqn:8.51}) for infinitesimal non-crossing cumulants now gives
\begin{equation}   \label{eqn:8.74}
\kappa_n ' (X_{i_1}, \ldots , X_{i_n}) = 
n \cdot \kappa_n (X_{i_1}, \ldots , X_{i_n}),
 \ \ \forall \,
n \geq 1, \ \forall \, 1 \leq i_1, \ldots ,i_n \leq k.
\end{equation}
In the particular case when $\mu$ is such that every $X_i$ is 
standard semicircular in $( \cA , \mu )$, it thus follows that 
every $X_i$ becomes a standard infinitesimal semicircular element in 
$( \cA , \mu , \muno ' )$, where we set the parameters from
Equation (\ref{eqn:5.71}) to be $\alpha_1 ' = 0$ and 
$\alpha_2 ' = 2$. On the other hand, if $\mu$ is defined such that 
every $X_i$ has a standard free Poisson distribution in 
$( \cA , \mu )$ (with $\kappa_n (X_i, \ldots , X_i) = 1$ for all 
$n \geq 1$), then the $X_i$ will become infinitesimal free Poisson 
elements in $( \cA , \mu , \muno ' )$, in the sense of Definition 
\ref{def:5.8} and where we take $\beta ' = 0$, $\gamma ' = 1$
in Equation (\ref{eqn:5.81}).
\end{example}

We now move to the situation described in Proposition 
\ref{prop:1.5}. Clearly, this is just an immediate 
consequence of Proposition \ref{prop:4.20}.

\begin{corollary}       \label{cor:8.8}
In the notations of Proposition \ref{prop:4.20}, suppose that 
$\cA_1, \ldots , \cA_k$ are unital subalgebras of $\cA$ which are 
freely independent with respect to $\varphi_t$ for every $t \in T$.
Then $\cA_1, \ldots , \cA_k$ are infinitesimally free in 
$( \cA , \varphi ,  \phiprim )$.
\end{corollary}

\begin{proof} Consider elements $a_1 \in \cA_{i_1}, \ldots ,
a_n \in \cA_{i_n}$ where the indices $i_1, \ldots , i_n$ are not 
all equal to each other. The freeness of $\cA_1, \ldots , \cA_k$ 
in $( \cA , \varphi_t )$ implies that 
$\kappa^{(t)}_n ( a_1, \ldots , a_n ) = 0$ for every $t \in T$.
The limit and derivative at $0$ of the function 
$t \mapsto \kappa^{(t)}_n ( a_1, \ldots , a_n )$ must therefore 
vanish, which means (by Proposition \ref{prop:4.20}) that 
$\kappa_n (a_1, \ldots , a_n) = \kappa_n ' (a_1, \ldots , a_n) =0$.
Hence condition (2) from Theorem \ref{thm:1.2} is satisfied, and
the conclusion follows.
\end{proof}

\begin{example}    \label{ex:8.9}
Consider again the situation where $\cA$ is the $*$-algebra
$\bC \langle X_1, \ldots , X_k \rangle$, as in Example \ref{ex:8.7},
and where $\mu : \cA \to \bC$ is a positive definite functional 
with $\mu (1) = 1$. Let $( \kappa_n )_{n \geq 1}$ be the 
non-crossing cumulant functionals of $\mu$, and let 
$\{ \kappaA_{\pi} \mid \pi \in \cup_{n=1}^{\infty} NC(n) \}$ be 
the extended family of multilinear functionals from Remark 
\ref{rem:3.14}.

For every $t > 0$, let $\mu_t : \cA \to \bC$ be the linear 
functional defined by putting $\mu_t (1) = 1$ and 
\begin{equation}   \label{eqn:8.91}
\mu_t ( X_{i_1}, \ldots , X_{i_n} ) = \sum_{\pi \in NC(n)} \,
(t+1)^{ | \pi | } \cdot \kappa_{\pi} ( X_{i_1}, \ldots , X_{i_n} ),
\end{equation}
for all $n \geq 1$ and $1 \leq i_1, \ldots , i_n \leq k$.
As is easily seen, $\mu_t$ is uniquely determined by the fact 
that its non-crossing cumulant functionals 
$( \, \kk_n^{(t)} \, )_{n \geq 1}$ satisfy
\begin{equation}   \label{eqn:8.92}
\kappa^{(t)}_n ( X_{i_1}, \ldots , X_{i_n} ) =
(t+1) \cdot \kappa_n ( X_{i_1}, \ldots , X_{i_n} ),
\ \ \forall \, n \geq 1, \ 1 \leq i_1, \ldots , i_n \leq k.
\end{equation}
Due to this fact, $\mu_t$ is called the ``$(t+1)$-th convolution
power of $\mu$'' with respect to the operation $\boxplus$ of 
free additive convolution -- see pp. 231-233 of \cite{NS2006} 
for details.

From (\ref{eqn:8.91}) it is clear that the family 
$\{ \mu_t \mid t > 0 \}$ has infinitesimal limit $( \mu , \mu ' )$ 
at $t = 0$, where $\mu$ is the functional we started with, while 
$\mu '$ is defined by putting $\mu ' (1) = 0$ and
\begin{equation}   \label{eqn:8.93}
\mu ' ( X_{i_1} \cdots  X_{i_n} ) = \sum_{\pi \in NC(n)} \,
| \pi | \cdot \kappa_{\pi} ( X_{i_1}, \ldots , X_{i_n} ),
\ \ \forall \, n \geq 1, \ 1 \leq i_1, \ldots , i_n \leq k.
\end{equation}
Note also that by using Equation (\ref{eqn:8.92}) and by invoking 
Proposition \ref{prop:4.20} we get 
\begin{equation}   \label{eqn:8.94}
\kappa_n ' ( X_{i_1}, \ldots , X_{i_n} ) = 
\kappa_n ( X_{i_1}, \ldots , X_{i_n} ),
\ \ \forall \, n \geq 1, \ 1 \leq i_1, \ldots , i_n \leq k.
\end{equation}

Now let $\cA_1, \ldots , \cA_k$ be the unital $*$-subalgebras of 
$\cA$ that were also considered in Example \ref{ex:8.7}, 
$\cA_i = \mbox{span} \{ X_i^n \mid n \geq 0 \}$ for 
$1 \leq i \leq k$. Suppose that $\cA_1, \ldots , \cA_k$ are 
free in $( \cA , \mu )$. Then they are free in $( \cA , \mu_t )$ 
for every $t > 0$; this follows from Equation (\ref{eqn:8.92}) 
and the description of freeness in terms of non-crossing 
cumulants (cf. Theorem 11.20 in \cite{NS2006}), where we take 
into account that $\cA_i$ is the unital algebra generated by 
$X_i$. Hence this is a situation where Corollary \ref{cor:8.8} 
applies, and we conclude that $\cA_1, \ldots , \cA_k$
are infinitesimally free in $( \cA , \mu, \mu ' )$.

Note also that if $X_i$ has a standard semicircular distribution 
in $( \cA , \mu )$, then Equation (\ref{eqn:8.94}) implies that 
$X_i$ becomes an infinitesimal semicircular element in 
$( \cA , \mu, \mu ' )$, where the parameters $\alpha_1 ', \alpha_2 '$
from Remark \ref{rem:5.7} are taken to be 
$\alpha_1 ' = 0$, $\alpha_2 ' = 1$. Likewise, if $X_i$ is a 
standard free Poisson in $( \cA , \mu )$, then Equation 
(\ref{eqn:8.94}) implies that $X_i$ becomes an infinitesimal 
free Poisson element in $( \cA , \mu, \mu ' )$, where the parameters 
$\beta ' , \gamma '$ from Definition \ref{def:5.8} are taken to be 
$\beta ' =1$, $\gamma ' = 0$. 
\end{example}

$\ $

{\bf\large Acknowledgements.} 
We are grateful to Teodor Banica for very useful discussions 
during the incipient stage of this work.
The second-named author would like to thank Mireille Capitaine 
and Muriel Casalis for inviting him to visit the Mathematics 
Institute of Universit\'e Paul Sabatier in Toulouse, where 
this work was started.

$\ $

$\ $

$\ $

Maxime F\'evrier

Institut de Math\'ematiques de Toulouse, Equipe de Statistique
et Probabilit\'es,

Universit\'e Paul Sabatier, 31062 Toulouse Cedex 09, France.

Email: fevrier@math.univ-toulouse.fr

$\ $

Alexandru Nica

Department of Pure Mathematics, University of Waterloo

Waterloo, Ontario N2L 3G1, Canada.

Email: anica@math.uwaterloo.ca

\end{document}